\PassOptionsToPackage{hyperfootnotes=false}{hyperref}
\documentclass[11pt]{article}

\usepackage{amsmath,amssymb,amsthm,mathtools,mathrsfs,bm}
\usepackage{graphicx,caption,booktabs,longtable,array,authblk}

\usepackage{tikz,float,needspace,microtype,placeins}
\usepackage[top=2.5cm,bottom=2.5cm,left=3cm,right=3cm]{geometry}
\usepackage[numbers,sort&compress]{natbib}
\setcitestyle{square,numbers,comma}
\usepackage{hyperref}
\hypersetup{
  pdftitle={Minimum-Norm Interpolation in Weighted Variation Spaces: Direct and Inverse Problems},
  pdfauthor={Weizhao Li, Lei Shi, Hanfei Zhou}
}

\allowdisplaybreaks[4]
\numberwithin{equation}{section}
\theoremstyle{plain}
\newtheorem{definition}{Definition}[section]
\newtheorem{theorem}{Theorem}[section]
\newtheorem{lemma}{Lemma}[section]
\newtheorem{proposition}{Proposition}[section]
\newtheorem{corollary}{Corollary}[section]
\newtheorem{example}{Example}[section]

\newtheorem{remark}[definition]{Remark}

\newcommand{\R}{\mathbb R}
\newcommand{\Sph}{\mathbb S}
\newcommand{\Mspace}{\mathcal M}
\newcommand{\Vspace}{\mathcal V}
\newcommand{\TV}{\mathrm{TV}}
\newcommand{\relu}{\operatorname{ReLU}}

\newcommand{\dist}{\operatorname{dist}}
\newcommand{\ran}{\operatorname{ran}}
\newcommand{\spann}{\operatorname{span}}
\newcommand{\argminn}{\operatorname*{argmin}}
\newcommand{\norm}[2][]{\left\lVert #2\right\rVert_{#1}}
\newcommand{\abs}[1]{\left\lvert #1\right\rvert}
\newcommand{\ip}[2]{\left\langle #1,#2\right\rangle}
\newcommand{\dd}[1]{\,\mathrm{d}#1}
\newcommand{\vect}[1]{\bm{#1}}
\newcommand{\matr}[1]{\bm{#1}}
\newcommand{\vx}{\vect{x}}
\newcommand{\vy}{\vect{y}}
\newcommand{\vu}{\vect{u}}
\newcommand{\vz}{\vect{z}}
\newcommand{\vxi}{\vect{\xi}}
\newcommand{\vtheta}{\vect{\theta}}
\newcommand{\vc}{\vect{c}}
\newcommand{\eps}{\varepsilon}

\newcommand{\Bd}{\mathbb{B}^d}
\newcommand{\Param}{\Theta}
\newcommand{\Vsp}{\Vspace_s}
\newcommand{\Vunw}{\Vspace_s^{0}(\Bd)}
\newcommand{\Msp}{\Mspace(\Param)}
\newcommand{\Ltwo}{L_2(\Bd)}
\newcommand{\Polys}{\Pi_s}

\newcommand{\rawatom}[2]{\phi_{#1,#2}}
\newcommand{\atom}[2]{\psi_{#1,#2}}
\newcommand{\paramatom}[1]{\psi_{#1}}
\newcommand{\wstar}[1]{w_s(#1)}
\newcommand{\synthop}{\mathscr T_s}
\newcommand{\synth}[1]{\synthop #1}
\newcommand{\vnorm}[1]{\norm[\Vsp]{#1}}
\newcommand{\unwnorm}[1]{\norm[\Vunw]{#1}}
\newcommand{\ltnorm}[1]{\norm[\Ltwo]{#1}}
\newcommand{\tvnorm}[1]{\norm[\TV]{#1}}

\newcommand{\netclass}[1]{\Sigma_{#1}}
\newcommand{\bestnerr}[1]{E_n(#1)_{\Ltwo}}
\newcommand{\weightparam}{\mathscr P_s}

\newcommand{\filldist}[1]{h_{#1,\Bd}}
\newcommand{\interpop}[1]{\mathscr I_{#1}}
\newcommand{\interp}[2]{\interpop{#1}#2}

\newcommand{\ixjf}{\interp{X_j}{f}}
\newcommand{\ixhf}{\interp{X_h}{f}}
\newcommand{\minset}[2]{\mathfrak M_{#1}(#2)}

\newcommand{\bernconst}[2]{\mathfrak B_{#1}(#2)}

\newcommand{\dict}[1]{\Lambda_{#1}}
\newcommand{\jointdict}[1]{\widehat{\Lambda}_{#1}}
\newcommand{\dictsynthop}[1]{\mathscr S_{#1}}
\newcommand{\dictsynth}[2]{\dictsynthop{#1}#2}
\newcommand{\wgram}[1]{\matr{G}_{#1}}
\newcommand{\mingrameig}[1]{\lambda_{\min}\!\left(\wgram{#1}\right)}
\newcommand{\privrad}[1]{\rho_{#1}}

\newcommand{\samplingop}[1]{\mathscr Q_{#1}}
\newcommand{\sampling}[1]{\samplingop{#1}}

\title{Minimum-Norm Interpolation in Weighted Variation Spaces: Direct and Inverse Problems\thanks{The work of Lei Shi is partially supported by the National Natural Science Foundation of China (Grant No.~12571099).
All authors contributed equally to this work and are listed in alphabetical order by surname.}}
\author[1]{Weizhao Li}
\author[1,2]{Lei Shi\thanks{Corresponding author. Email: \texttt{leishi@fudan.edu.cn}.}}
\author[1]{Hanfei Zhou}
\affil[1]{School of Mathematical Sciences, \linebreak
Fudan University, Shanghai, 200433, China}
\affil[2]{Shanghai Key Laboratory for Contemporary Applied Mathematics, \linebreak
Fudan University, Shanghai, 200433, China \linebreak
Email: weizhaoli24@m.fudan.edu.cn, leishi@fudan.edu.cn, zhouhf23@m.fudan.edu.cn}
\date{}

\begin{document}
\maketitle

\begin{abstract}
We study minimum-norm interpolation in weighted variation spaces generated by ReLU$^s$ ridge functions on the unit ball, where $s$ is a positive integer.
For a target \(f\) in the weighted variation space and a node set with sufficiently small fill distance $h$, every minimum-norm interpolant satisfies the direct bound
\[
 \norm[\Ltwo]{f-\ixhf}\lesssim_{d,s}h^{(2s+1)/4}\vnorm{f}.
\]
Under a boundary separation condition on the nodes, we construct a matching lower bound for the worst-case interpolation error.
For the inverse problem, we introduce dual certificates and derive an estimate analogous to a Bernstein inequality.
Under suitable upper bounds on dictionary size and spectral lower bounds for the atomic Gram matrices, sufficiently fast \(L_2\) convergence of minimum-norm interpolants along geometrically decreasing fill distances implies that the target belongs to \(\Vsp\).
\leavevmode\phantomsection\label{not:network-error}As a byproduct, the best \(L_2\) approximation error \(\bestnerr{f}\) for shallow ReLU$^s$ networks of width at most \(n\) satisfies
\[
 \bestnerr{f}\lesssim_{d,s}n^{-1/2-(2s+1)/(2d)}\vnorm{f},
\]
with the weighted atomic cost of an approximating network bounded by a constant multiple of \(\vnorm{f}\), independently of \(n\).
\end{abstract}


\section{Introduction}\label{sec:introduction}

Although neural networks have achieved considerable success in image recognition and natural language processing \citep{he2016deep,vaswani2017attention}, understanding their generalization properties remains a central problem.
In particular, overparameterized networks, whose number of parameters exceeds the number of training samples, can fit the observed data almost exactly and still have a small prediction error \citep{zhang2017rethinking,belkin2019reconciling,neyshabur2019role,arora2019fine,nakkiran2020double}.
One influential explanation is that training dynamics can favor solutions with small norm, even in the absence of an explicit penalty; this selection is known as implicit regularization \citep{bartlett2021deep,hastie2022surprises,boursier2022gradient}.
This phenomenon motivates minimum-norm interpolation as a model for selecting a predictor from the functions that fit the data.

Let \(\Omega\subset\R^d\) be a bounded domain, where \(d\ge1\), and let \(\mathcal F\subset C(\Omega)\cap L_2(\Omega)\) be a normed function space.
\leavevmode\phantomsection\label{not:lp-spaces}We use the standard Lebesgue spaces \(L_p(\Omega)\), \(1\le p\le\infty\), with norms \(\norm[L_p(\Omega)]{\cdot}\).
Given point observations \(f(\vx_i)\) of a target \(f:\Omega\to\R\) in \(L_2(\Omega)\) at distinct nodes \(X=\{\vx_1,\ldots,\vx_m\}\subset\Omega\), we consider the problem
\begin{equation}\label{eq:min-interpolation}
 \min_{g\in\mathcal F}\bigl\{\norm[\mathcal F]{g}:g(\vx_i)=f(\vx_i),\ 1\le i\le m\bigr\}.
\end{equation}
We say that the node set $X$ is feasible if problem \eqref{eq:min-interpolation} has a solution.
In \eqref{eq:min-interpolation}, the choice of model space \(\mathcal F\) and its norm determines the cost assigned to an interpolating function.
We will specify this space below.

We focus on two complementary questions relating approximation rates to the regularity or representation properties of the target:
\begin{itemize}
 \item \emph{Direct problem:} Given a function class for the target, at what rate does the approximation error decrease as the interpolation nodes become dense or the network width increases?
 \item \emph{Inverse problem:} Given a rate of convergence of the approximants, what can be inferred about the function space to which the target belongs?
\end{itemize}
For minimum-norm interpolation, we study these questions through the error \(\norm[L_2(\Omega)]{f-g_X}\), where \(g_X\) solves \eqref{eq:min-interpolation}.
The geometry of the nodes is measured by the fill distance and, when \(m\ge2\), the separation distance:
\begin{equation}\label{eq:fill-distance}
 h=h_{X,\Omega}:=\sup_{\vx\in\Omega}\min_{\vy\in X}\abs{\vx-\vy},\ \  q_X:=\frac12\min_{\substack{\vx,\vy\in X\\\vx\ne\vy}}\abs{\vx-\vy}.
\end{equation}
A family of node sets is called quasi-uniform if \(h_{X,\Omega}/q_X\) is bounded uniformly over the family.
We use the fill distance \(h\) to quantify convergence rates in the direct and inverse problems.

Classical estimates for scattered data interpolation relate these rates to Sobolev regularity and the geometry of the nodes.
\leavevmode\phantomsection\label{not:sobolev-spaces}We write \(H^\beta(\Omega)\) for the \(L_2\) Sobolev space of order \(\beta\), including fractional orders \citep{adamsfournier2003sobolev}.
Bounds for functions that vanish at scattered nodes control their lower order Sobolev norms in terms of the fill distance and a higher order norm \citep{narcowich2005zeros}.
For example, on a compact Lipschitz region \(\Omega\), consider a radial basis function kernel whose Fourier transform is comparable to \((1+|\omega|^2)^{-\tau}\), where \(\tau>d/2\).
For every \(f\in H^\beta(\Omega)\) with \(d/2<\beta\le\tau\), the error bound
\begin{equation*}
 \norm[L_2(\Omega)]{f-g_X}=O(h^\beta)
\end{equation*}
holds for kernel interpolants on quasi-uniform node sets \citep{narcowich2006sobolev,wenzel2025sharp}.
For the inverse problem, let \(f\in C(\Omega)\).
A recent inverse theorem \citep{wenzel2025sharp} gives
\begin{equation*}
 \norm[L_2(\Omega)]{f-g_X}=O(h^\beta)
 \quad\Longrightarrow\quad
 \begin{cases}
  f\in H^{\beta'}(\Omega), & 0<\beta'<\beta\le\tau,\\
  f\in H^\tau(\Omega), & \beta>\tau.
 \end{cases}
\end{equation*}
Here the error bound is uniform over all sufficiently dense node sets satisfying \(h/q_X\le44\).
These results are established for classical kernel interpolation and do not address minimum-norm interpolation by neural networks.

For targets in the Barron space, shallow ReLU networks of width \(n\) attain an \(L_2\) approximation error \(O(n^{-1/2})\); conversely, continuous pointwise limits of networks with uniformly bounded path norm belong to this space \citep{emawu2022barron}.
For linearized shallow \(\relu^s\) networks on \(\Sph^d\), discrete least squares with prescribed hidden parameters attains an \(L_2\) error \(O(n^{-r/d})\) for targets in \(H^r(\Sph^d)\) under suitable assumptions \citep{liumaoxu2026deterministic}.
The same work also gives an inverse result relating best approximation rates to Sobolev regularity.
These two works address direct and inverse network approximation, but do not impose interpolation constraints.
Minimum complexity interpolation with \(L_p\) coefficient costs, \(p>1\), is studied in \citep{celentano2021minimum}, where finite random features interpolants are shown to approximate the corresponding integral interpolant under suitable assumptions.
That work considers the limit as the number of features grows for fixed data, rather than the direct and inverse problems associated with increasingly dense interpolation nodes.

In this paper, we study direct and inverse problems for minimum-norm interpolation by shallow \(\relu^s\) networks, allowing both the hidden parameters and the output coefficients to vary.
We take \(\mathcal F\) in \eqref{eq:min-interpolation} to be a weighted variation space \(\Vsp\), following the construction for ReLU networks in \citep{devore2025weighted}.
This space consists of functions admitting integral network representations with finite weighted total variation, and its norm is the infimum of this cost over all such representations.
The full definition is given in Section~\ref{sec:preliminaries}.
\leavevmode\phantomsection\label{not:domain}We work on the unit ball \(\Bd:=\{\vx\in\R^d:\abs{\vx}<1\}\) and use the activation \(\relu^s(z)=z_+^s\) for a fixed integer \(s\ge1\).
\leavevmode\phantomsection\label{not:network-class}Write
\begin{equation}\label{eq:network-class}
 \netclass{n}:=\left\{\sum_{j=1}^n a_j(\vxi_j\cdot\vx-t_j)_+^s:
 a_j\in\R,\ \vxi_j\in\Sph^{d-1},\ -1\le t_j<1\right\}
\end{equation}
for the shallow networks of width at most \(n\).
This simple architecture is a useful model for studying network representation costs and approximation rates \citep{emawu2022barron,siegelxu2023variation,siegelxu2024sharp,liliushi2026shallow,yangzhou2025optimal}.
\leavevmode\phantomsection\label{not:weight}\label{not:atomic-cost}We assign each such representation the weighted coefficient cost
\begin{equation}\label{eq:weighted-atomic-cost}
 \sum_{j=1}^n\wstar{t_j}\abs{a_j},
 \qquad \wstar{t}:=(1-t)^{(d+2s)/4}.
\end{equation}
The weighted cost in \eqref{eq:weighted-atomic-cost} accounts for the support of each ridge function: as \(t\uparrow1\), its supporting spherical cap shrinks and the weight \(\wstar{t}\) decreases, assigning a smaller cost to functions supported on small caps.

Our main results are as follows.
{\renewcommand{\labelenumi}{\textup{(\arabic{enumi})}}%
\setlength{\leftmargini}{3em}%
\begin{enumerate}
 \item For every integer \(s\ge1\) and every node set \(X_h\) with sufficiently small fill distance \(h\), each minimum-norm interpolant of \(f\in\Vsp\) satisfies
 \begin{equation*}
  \norm[\Ltwo]{f-\ixhf}\lesssim_{d,s}h^{(2s+1)/4}\vnorm{f}.
 \end{equation*}
 Under the boundary separation condition \eqref{eq:interior-margin}, which requires each node to stay at distance at least \(\kappa h\) from the boundary for some fixed \(\kappa>0\), we construct functions in the unit ball of \(\Vsp\) that vanish at all the nodes and have \(L_2\) norm of this order.
 This gives a matching lower bound and establishes the sharpness of the interpolation estimate.
 \item We introduce dual certificates to obtain an estimate analogous to a Bernstein inequality and prove inverse theorems.
 Specifically, let \(X_j\) be increasingly dense node sets with fill distances \(h_j\) decreasing geometrically to zero, and let \(g_j\) be the corresponding minimum-norm network interpolants of \(f\).
 We prove that if \(g_j\) and \(g_{j+1}\) together use \(O(h_{j+1}^{-\nu})\) distinct hidden parameters and the \(L_2\) Gram matrix of their normalized atoms has smallest eigenvalue \(\gtrsim h_{j+1}^{2\eta}\), where \(\nu,\eta\ge0\), then
 \begin{equation*}
  \ltnorm{f-g_j}=O(h_j^\beta),\qquad \beta>\frac\nu2+\eta
  \quad\Longrightarrow\quad f\in\Vsp.
 \end{equation*}
 \item We prove the estimate
 \begin{equation*}
    \inf_{g\in\netclass{n}}\ltnorm{f-g}\lesssim_{d,s}n^{-1/2-(2s+1)/(2d)}\vnorm{f},
 \end{equation*}
 with a uniformly bounded weighted sum of the absolute values of the output coefficients.
 For \(s=1\), the bound \(O(n^{-1/2-3/(2d)})\) on the unit ball of the weighted variation space was obtained in \citep{devore2025weighted}; our result extends it to every integer \(s\ge1\).
\end{enumerate}
}

The remainder of the paper is organized as follows.
Section~\ref{sec:preliminaries} introduces weighted variation spaces, minimum-norm interpolation, and finite dictionary certificates.
Section~\ref{sec:main-results} presents the direct theorem, the optimality result, and inverse theorems under different stability conditions.
Section~\ref{sec:discussion-conclusion} summarizes the main conclusions and discusses open questions.
Auxiliary results and the proofs of the direct and inverse theorems are collected in the appendices.
For the reader's convenience, Appendix~\ref{app:notation} provides a reference table of the principal notation used in the paper.
\section{Preliminaries}\label{sec:preliminaries}
We first introduce weighted variation spaces and formulate the minimum-norm interpolation problem in this setting.
We then present the dual certificates used in the inverse theorems.

\subsection{Weighted variation spaces}
We extend the weighted variation space for ReLU networks introduced in \citep{devore2025weighted} to the activation \(\relu^s\), for any integer \(s\ge1\).
Further properties and proofs are given in Appendix~\ref{app:weighted-network-auxiliary-results}.

\leavevmode\phantomsection\label{not:parameter-domain}With \(\Bd\), \(s\), and \(\wstar{t}\) as in Section~\ref{sec:introduction}, let \(\Param:=\Sph^{d-1}\times[-1,1]\).
\leavevmode\phantomsection\label{not:raw-atom}For \(\vtheta=(\vxi,t)\in\Param\), define the \(\relu^s\) ridge atom by
\[
 \rawatom{\vxi}{t}(\vx):=(\vxi\cdot\vx-t)_+^s.
\]
For \(-1\le t<1\), normalize this atom using the weight in \eqref{eq:weighted-atomic-cost}:
\begin{equation}\label{eq:normalized-cap-atom}
 \atom{\vxi}{t}(\vx) :=\frac{\rawatom{\vxi}{t}(\vx)}{\wstar{t}} =\frac{(\vxi\cdot\vx-t)_+^s}{(1-t)^{(d+2s)/4}},
\end{equation}
and set \(\atom{\vxi}{1}=0\) at the endpoint \(t=1\).
\leavevmode\phantomsection\label{not:radon-measures}Let \(\Msp\) denote the space of finite signed Radon measures on \(\Param\), equipped with the total variation norm \(\tvnorm{\cdot}\) \citep{folland1999real}.
\leavevmode\phantomsection\label{not:synthesis}Define the synthesis operator by
\begin{equation*}
 \synthop:\Msp\longrightarrow\Ltwo, \qquad \synth{\mu} :=\int_{\Param}\paramatom{\vtheta}\dd{\mu}(\vtheta).
\end{equation*}

\leavevmode\phantomsection\label{not:weighted-space}The weighted variation space and its norm are defined by
\begin{equation*}
 \Vsp:=\Vspace_{w_s}(\Bd):=\ran(\synthop),\qquad \vnorm{f}:=\inf_{\synth{\mu}=f}\tvnorm{\mu}.
\end{equation*}
\leavevmode\phantomsection\label{not:network-coefficients}A network in \eqref{eq:network-class} can be written as
\[
 g=\sum_{j=1}^n c_j\atom{\vxi_j}{t_j}, \qquad c_j=\wstar{t_j}a_j,
\]
and its norm satisfies
\begin{equation*}
 \vnorm{g} \le\sum_{j=1}^n\abs{c_j} =\sum_{j=1}^n\wstar{t_j}\abs{a_j},
\end{equation*}
where the last sum is the coefficient cost in \eqref{eq:weighted-atomic-cost}.
The norm \(\vnorm{g}\) is the least cost over all measure representations of the function, so it provides the objective in \eqref{eq:min-interpolation} independently of a particular network representation.

Every element of \(\Vsp\) has a unique continuous representative on \(\Bd\), and point evaluation at each interior point is a continuous linear functional on \(\Vsp\); see Appendix~\ref{app:basic-properties}.
We use this representative when imposing interpolation conditions.
The following closure property will be used in the inverse theorems.

\begin{proposition}
\label{thm:weak-star-closure}
Let \(g_n\in\Vsp\) satisfy
\begin{equation*}
 g_n\to f\quad\text{in }\Ltwo, \qquad \sup_n\vnorm{g_n}<\infty.
\end{equation*}
Then \(f\in\Vsp\), and
\begin{equation*}
 \vnorm{f} \le\liminf_{n\to\infty}\vnorm{g_n}.
\end{equation*}
\end{proposition}

Together with Theorem~\ref{thm:weighted-network-approximation}, this shows that \(\Vsp\) consists precisely of the \(\Ltwo\) limits of shallow networks whose weighted coefficient costs remain uniformly bounded.
Thus, to identify the limit of a sequence of interpolants as an element of \(\Vsp\), it suffices to bound their variation norms.
The dual certificates introduced below provide such a bound in Section~\ref{sec:inverse}.

\subsection{Minimum-norm interpolation problem}\label{sec:interpolation-duality}
\leavevmode\phantomsection\label{not:node-data}We now specialize \eqref{eq:min-interpolation} to \(\mathcal F=\Vsp\) and \(X\subset\Bd\).
Data \(\vy=(y_1,\ldots,y_m)\) are called feasible if some \(g\in\Vsp\) satisfies \(g(\vx_i)=y_i\) for \(1\le i\le m\).
\leavevmode\phantomsection\label{not:minimizer-set}For data generated by a function \(f\), define the minimizer set by
\begin{equation*}
 \mathfrak M_X(f):= \argminn_{g\in\Vsp}\bigl\{\vnorm{g}:g|_X=f|_X\bigr\},
\end{equation*}
\leavevmode\phantomsection\label{not:interpolant}and use \(\interp{X}{f}\in\mathfrak M_X(f)\) to denote an arbitrary minimum-norm interpolant.
We next state the corresponding representer theorem; see \citep[Theorems~3.3 and~4.2]{brediescarioni2020sparsity}.

\begin{proposition}
\label{prop:representer}
For feasible data, problem \eqref{eq:min-interpolation} with \(\mathcal F=\Vsp\) has a solution of the form
\begin{equation*}
 g(\vx)=\sum_{j=1}^{n}a_j(\vxi_j\cdot\vx-t_j)_+^s, \qquad n\le m,
\end{equation*}
where the representation attains the minimum weighted atomic cost.
\end{proposition}
The representer theorem gives a network realization of the direct estimate in Section~\ref{sec:direct}.
For the inverse problem, a representation attaining the minimum coefficient cost allows us to express the norm of the interpolant through a finite \(\ell^1\) minimization problem and its dual.
\subsection{Finite dictionary certificates}\label{sec:finite-certificates}
In this subsection, we introduce dual certificates, a key tool in the proofs of the inverse theorems.
Related dual certificates are used in super-resolution \citep{fernandezgranda2016superresolution} and in the analysis of sparsity for shallow ReLU networks with total variation regularization \citep{delgrandebrune2026certificate}.

\leavevmode\phantomsection\label{not:dictionary-matrices}For a finite parameter set \(\Lambda=\{\vtheta_1,\ldots,\vtheta_M\}\subset\Param\) whose corresponding atoms are nonzero, define
\begin{equation*}
 \dictsynth{\Lambda}{\vc}:=\sum_{k=1}^M c_k\paramatom{\vtheta_k}, \qquad
 \wgram{\Lambda}:=\left(\ip{\paramatom{\vtheta_k}}{\paramatom{\vtheta_\ell}}_{\Ltwo}\right)_{k,\ell=1}^M,
\end{equation*}
and let \(E=\spann\{\paramatom{\vtheta}:\vtheta\in\Lambda\}\).
Suppose that the representation \(g=\dictsynth{\Lambda}{\vc}\) attains the \(\Vsp\) norm, so that
\begin{equation*}
 \vnorm{g}=\inf_{\synth{\mu}=g}\tvnorm{\mu}=\min_{\vect b\in\R^M}\left\{\norm[\ell^1]{\vect b}:\dictsynth{\Lambda}{\vect b}=g\right\}=\norm[\ell^1]{\vc}.
\end{equation*}
Duality for this finite \(\ell^1\) minimization problem gives
\begin{equation}\label{eq:finite-primal-problem}
 \vnorm{g}=\max_{q\in E}\left\{\ip{q}{g}_{\Ltwo}:\left|\ip{q}{\paramatom{\vtheta}}_{\Ltwo}\right|\le1\quad(\vtheta\in\Lambda)\right\}.
\end{equation}
A feasible \(q\) in \eqref{eq:finite-primal-problem} is called a dual certificate for this representation problem if \(\ip{q}{g}_{\Ltwo}=\vnorm{g}\).
For such a certificate, if \(u=\dictsynth{\Lambda}{\vect b}\) also satisfies \(\norm[\ell^1]{\vect b}=\vnorm{u}\), then \(\ip{q}{u}_{\Ltwo}\le\vnorm{u}\), and hence
\begin{equation}\label{eq:certificate-increment}
 \vnorm{g}-\vnorm{u}\le\ip{q}{g-u}_{\Ltwo}\le\ltnorm{q}\,\ltnorm{g-u}.
\end{equation}
This bound on the variation norm increment plays a role analogous to that of Bernstein inequalities, which control stronger Sobolev norms by weaker norms on suitable approximation spaces.
Such inequalities are used in inverse approximation results and convergence analyses of least squares methods \citep{wenzel2025sharp,liumaoxu2026deterministic,mao2026convergence}.
The factor \(\ltnorm{q}\) in \eqref{eq:certificate-increment} motivates choosing a certificate of minimum \(\Ltwo\) norm.
Therefore, we define
\begin{equation}\label{eq:dictionary-certificates}
 \begin{split}
 \mathcal C_\Lambda(g)&:=\left\{q\in E:\left|\ip{q}{\paramatom{\vtheta}}_{\Ltwo}\right|\le1\ (\vtheta\in\Lambda),\quad\ip{q}{g}_{\Ltwo}=\vnorm{g}\right\},\\
 \Gamma_\Lambda(g)&:=\min_{q\in\mathcal C_\Lambda(g)}\ltnorm{q}.
 \end{split}
\end{equation}
For \(\Lambda=\varnothing\) and \(g=0\), set \(\mathcal C_\Lambda(g)=\{0\}\) and \(\Gamma_\Lambda(g)=0\).
The following proposition characterizes \(\Gamma_\Lambda(g)\) and bounds it in terms of the Gram matrix.

\begin{proposition}\label{prop:finite-certificates}
Under the preceding assumptions, \(\mathcal C_\Lambda(g)\) is nonempty, compact, and convex, and its minimum-norm element is unique.
Let \(G=\wgram{\Lambda}\), and denote its Moore--Penrose pseudoinverse by \(G^\dagger\).
Then
\begin{equation}\label{eq:certificate-quadratic-program}
 \Gamma_\Lambda(g)^2=\min\left\{\vect a^{\mathsf T}G^\dagger\vect a:
 \begin{array}{l}
 \vect a\in\ran G,\quad\norm[\ell^\infty]{\vect a}\le1,\quad a_k=\operatorname{sgn}(c_k)\quad\text{if }c_k\ne0
 \end{array}\right\}.
\end{equation}
If \(G\) is positive definite, then
\begin{equation}\label{eq:certificate-spectral-bound}
 \Gamma_\Lambda(g)\le\sqrt{\frac{\abs{\Lambda}}{\mingrameig{\Lambda}}}.
\end{equation}
\end{proposition}
Thus, the minimum certificate norm is given by a convex quadratic program, which remains valid when the Gram matrix is singular.
The proof of Proposition~\ref{prop:finite-certificates} is given in Appendix~\ref{app:inverse-proofs}.

\section{Main Results}\label{sec:main-results}
We present the main results in this section and defer their proofs to the appendices.
Section~\ref{sec:direct} gives the direct interpolation estimate and a matching lower bound under the boundary separation condition.
Section~\ref{sec:inverse} establishes inverse theorems using dual certificates and lower spectral bounds for the weighted Gram matrix.

\subsection{The Direct Theorem and Its Sharpness}\label{sec:direct}
For a finite node set \(X_h\subset\Bd\), write \(h=\filldist{X_h}\), with the fill distance defined in \eqref{eq:fill-distance}.
The following theorem gives an \(L_2\) error bound for minimum-norm interpolants.
Its proof is deferred to Appendix \ref{app:direct-proof}.
\begin{theorem}\label{thm:direct}
Let \(X_h\subset\Bd\) have sufficiently small fill distance \(h\).
If \(f\in\Vsp\) and \(\ixhf\in\mathfrak M_{X_h}(f)\), then
\begin{equation*}
 \norm[\Ltwo]{f-\ixhf} \lesssim_{d,s}h^{(2s+1)/4}\vnorm{f}.
\end{equation*}
\end{theorem}

\begin{remark}
Minimality gives \(\vnorm{\ixhf}\le\vnorm{f}\), which is used in the proof to control the error away from the nodes.
In particular, \(\ixhf\to f\) in \(\Ltwo\) as \(h\to0\).
\end{remark}

\begin{remark}
Applying the sampling inequality in the proof to functions that vanish at the nodes shows that if a network \(g_h\in\Vsp\) satisfies
\begin{equation*}
 g_h|_{X_h}=f|_{X_h},\qquad \vnorm{g_h}\le C_0\vnorm{f},
\end{equation*}
where \(C_0\) is independent of \(h\), then
\begin{equation*}
 \norm[\Ltwo]{f-g_h} \lesssim_{d,s}(1+C_0)h^{(2s+1)/4}\vnorm{f}.
\end{equation*}
Thus, the same convergence rate holds whenever the weighted variation norm of the interpolant is bounded by a fixed multiple of that of the target function.
\end{remark}


We now consider a lower bound for the interpolation error.
For this purpose, define \leavevmode\phantomsection\label{not:sampling-radius}the quantity
\begin{equation*}
 \mathcal R_s(X_h):=\sup\bigl\{ \ltnorm{u}:\vnorm{u}\le1,\ u|_{X_h}=0 \bigr\}.
\end{equation*}
The interpolation property and minimality give
\begin{equation*}
 \mathcal R_s(X_h)\le\sup_{\vnorm{f}\le1}\norm[\Ltwo]{f-\ixhf}\le2\mathcal R_s(X_h).
\end{equation*}
The lower bound follows by considering targets that vanish on \(X_h\), whose minimum-norm interpolant is zero.
For the upper bound, \((f-\ixhf)/2\) vanishes on \(X_h\) and has variation norm at most one whenever \(\vnorm{f}\le1\).
Fix \(0<\kappa<1\).
We say that the node set satisfies the boundary separation condition if
\begin{equation}\label{eq:interior-margin}
 \abs{\vx_i}\le1-\kappa h,\qquad \vx_i\in X_h.
\end{equation}
This condition places every node at distance at least \(\kappa h\) from the boundary.
Our direct result has a matching lower bound in the following sense.
\begin{theorem}\label{thm:direct-lower}
Let \(X_h\subset\Bd\) have sufficiently small fill distance \(h\) and satisfy the boundary separation condition \eqref{eq:interior-margin}.
Then
\begin{equation*}
 \mathcal R_s(X_h)\asymp_{d,s,\kappa}h^{(2s+1)/4}.
\end{equation*}
\end{theorem}
This establishes the optimal order of the worst-case interpolation error over the unit ball of \(\Vsp\) for node sets satisfying \eqref{eq:interior-margin}.
The proof is given in Appendix~\ref{app:direct-proof}.

Finally, as a part of the direct problem, we also record the following network approximation result, proved in Appendix~\ref{app:weighted-network-auxiliary-results}.
For $s=1$, this recovers the approximation theorem in \citep{devore2025weighted}.
\begin{theorem}
\label{thm:weighted-network-approximation}
For every \(f\in\Vsp\) and \(n\ge1\),
\begin{equation*}
 \inf_{g\in\netclass{n}}\ltnorm{f-g} \lesssim_{d,s} n^{-1/2-(2s+1)/(2d)}\vnorm{f}.
\end{equation*}
Moreover, a network representation attaining this rate can be chosen so that its weighted atomic cost satisfies
\begin{equation*}
 \sum_{j=1}^n\wstar{t_j}\abs{a_j} \lesssim_{d,s}\vnorm{f}.
\end{equation*}
\end{theorem}

\subsection{The Inverse Problem}\label{sec:inverse}
We now study conditions under which the decay of the minimum-norm interpolation error implies that the target belongs to the weighted variation space.
The proofs are given in Appendices~\ref{app:inverse-proofs} and~\ref{app:private-patches}.

\leavevmode\phantomsection\label{not:node-sequence}Let \(\mathbf X=(X_j)_{j\ge1}\) be a sequence of finite node sets in \(\Bd\), and write \(h_j:=\filldist{X_j}\).
We consider scales that decrease geometrically, in the sense that there are constants \(0<\underline\rho\le\overline\rho<1\) such that
\begin{equation}\label{eq:geometric-scales}
 \underline\rho h_j\le h_{j+1}\le\overline\rho h_j, \qquad j\ge1.
\end{equation}
\begin{remark}
A sequence satisfying \eqref{eq:geometric-scales} can be constructed by taking \(\delta_j=4^{-j}\) and choosing \(X_j\subset\Bd\) to be maximal with respect to inclusion among sets whose distinct points are at least \(\delta_j\) apart.
Such sets are finite, and maximality gives \(h_j\le\delta_j\).
The midpoint of a closest pair of nodes lies in \(\Bd\) and is at distance at least \(\delta_j/2\) from every node, so \(h_j\ge\delta_j/2\).
Consequently,
\[
 \frac18h_j\le h_{j+1}\le\frac12h_j.
\]
The inverse results below require \eqref{eq:geometric-scales}, but do not require the node sets to be nested or quasi-uniform.
\end{remark}
\leavevmode\phantomsection\label{not:inverse-errors}Let \(f:\Bd\to\R\) belong to \(\Ltwo\), and fix its point values on \(\bigcup_{j\ge1}X_j\).
Assume that the data \(f|_{X_j}\) are feasible in \(\Vsp\), and abbreviate \(\minset{j}{f}:=\mathfrak M_{X_j}(f)\).
For \(g_j\in\minset{j}{f}\), write \(e_j=\norm[\Ltwo]{f-g_j}\) and \(m_j:=\vnorm{g_j}\).
\leavevmode\phantomsection\label{not:minimum-norm}
By Proposition \ref{thm:weak-star-closure}, when \(e_j\to0\), uniform boundedness of \((m_j)\) is sufficient to ensure that \(f\in\Vsp\).

We apply the dual certificates of Section~\ref{sec:finite-certificates} to successive interpolants.
By Proposition \ref{prop:representer}, a minimum-norm interpolant at each scale can be chosen as a finite network whose representation attains the minimum weighted atomic cost.
After merging repeated parameters and removing zero atoms, denote the parameter set by \(\dict{j}\) and write
\begin{equation}\label{eq:minimum-cost-representations}
 g_j=\sum_{\vtheta\in\dict{j}}c_{j,\vtheta}\paramatom{\vtheta}, \qquad m_j=\sum_{\vtheta\in\dict{j}}\abs{c_{j,\vtheta}}.
\end{equation}
\leavevmode\phantomsection\label{not:joint-dictionary}Set \(\jointdict{j}=\dict{j}\cup\dict{j+1}\) and \(E_j=\spann\{\paramatom{\vtheta}:\vtheta\in\jointdict{j}\}\).
Using \eqref{eq:dictionary-certificates}, set
\begin{equation}\label{eq:finite-certificate-set}
 \mathcal C_j:=\mathcal C_{\jointdict{j}}(g_{j+1}),\qquad
 \Gamma_j:=\Gamma_{\jointdict{j}}(g_{j+1}).
\end{equation}
The following theorem gives sufficient rates of decay for \(e_j\) in terms of the growth of \(\Gamma_j\).

\begin{theorem}
\label{thm:certificate-inverse}
Suppose that the scales satisfy \eqref{eq:geometric-scales}, and choose representations as in \eqref{eq:minimum-cost-representations}.
Assume that the certificate norms in \eqref{eq:finite-certificate-set} satisfy
\begin{equation*}
 \Gamma_j\le C_\Gamma h_{j+1}^{-\mu},\qquad \mu\ge0.
\end{equation*}
\leavevmode\phantomsection\label{not:inverse-rates}If, for some \(A>0\), either
\begin{equation}\label{eq:critical-log-rate}
 \begin{aligned}
 &e_j\le Ah_j^\beta,\qquad\beta>\mu,\\
 \text{or}\quad &e_j\le Ah_j^\mu(1+j)^{-1-\eps},\qquad\eps>0,
 \end{aligned}
\end{equation}
then \(f\in\Vsp\).
\end{theorem}
\begin{remark}
  Under \eqref{eq:geometric-scales}, the second error condition in \eqref{eq:critical-log-rate} is equivalent, up to constants for sufficiently large \(j\), to \(e_j\lesssim h_j^\mu|\log h_j|^{-1-\eps}\).
\end{remark}
\begin{remark}\label{cor:approximate-certificates}
The inverse conclusion also holds for approximate dual certificates.
Let \(q_j\in E_j\) and \(\tau_j\ge0\) satisfy
\begin{equation*}
 \abs{\ip{q_j}{\paramatom{\vtheta}}_{\Ltwo}}\le1\quad(\vtheta\in\jointdict{j}),
 \qquad \ip{q_j}{g_{j+1}}_{\Ltwo}\ge m_{j+1}-\tau_j.
\end{equation*}
If \(\sum_{j=1}^\infty\tau_j<\infty\) and \(\ltnorm{q_j}\le C_qh_{j+1}^{-\mu}\) for some \(C_q>0\) and \(\mu\ge0\), either error condition in \eqref{eq:critical-log-rate} implies \(f\in\Vsp\).
\end{remark}

The spectral estimate \eqref{eq:certificate-spectral-bound} leads to an inverse result stated directly in terms of the dictionary cardinality and the smallest eigenvalue of the weighted Gram matrix.
\begin{corollary}\label{thm:gram-inverse}
Suppose that the scales satisfy \eqref{eq:geometric-scales}, and choose representations as in \eqref{eq:minimum-cost-representations}.
\leavevmode\phantomsection\label{not:gram-exponents}Assume that there are \(C_M,c_G>0\) and \(\nu,\eta\ge0\) such that
\begin{equation}\label{eq:spectral-condition}
 \abs{\jointdict{j}}\le C_Mh_{j+1}^{-\nu},\qquad
 \mingrameig{\jointdict{j}}\ge c_Gh_{j+1}^{2\eta}.
\end{equation}
If
\begin{equation*}
 e_j\le Ah_j^\beta,\qquad\beta>\frac\nu2+\eta,
\end{equation*}
then \(f\in\Vsp\).
\end{corollary}

\begin{remark}
For node families satisfying \(\abs{X_j}\lesssim h_j^{-d}\), Proposition \ref{prop:representer} ensures that each interpolating network can be chosen with width no greater than the number of nodes.
One may therefore take \(\nu=d\) in the cardinality condition.
\end{remark}

Corollary~\ref{thm:gram-inverse} follows directly from Proposition~\ref{prop:finite-certificates} and Theorem~\ref{thm:certificate-inverse}.
The alternative proof in Appendix~\ref{app:inverse-proofs} controls the \(\Vsp\) norms of successive differences by their \(\Ltwo\) norms and shows that \((g_j)\) is Cauchy in \(\Vsp\).
A related argument for Sobolev kernel interpolation uses a spectral bound for the kernel matrix \citep[Theorem~5 and proof of Theorem~1]{wenzel2025sharp}.
To be precise, for a continuous kernel \(k\) with native space \(\mathcal H_k(\mathbb{B}^d)=H^\tau(\mathbb{B}^d)\) and equivalent norms, where \(\tau>d/2\), this bound is
\begin{equation*}
 A_X:=\bigl(k(\vx_i,\vx_\ell)\bigr)_{i,\ell=1}^m,
 \qquad \lambda_{\min}(A_X)\gtrsim q_X^{2\tau-d}.
\end{equation*}
Here the separation distance concerns the interpolation nodes, whereas \eqref{eq:spectral-condition} bounds the \(L_2\) Gram matrix of the selected network atoms.

For shallow networks with prescribed hidden parameters, inverse inequalities in Gaussian Sobolev norms depend on parameter separation and the bias range \citep{mao2026convergence}.
When the hidden parameters are also optimized, nearby atoms can nearly cancel, and atoms supported on shrinking boundary caps have small \(L_2\) norms.
Consequently, network width alone cannot control the \(\Vsp\) norm by the \(\Ltwo\) norm; see Appendix~\ref{app:inverse-proofs}.
The following geometric condition on the hidden parameters controls both effects and gives a spectral lower bound, leading to an inverse theorem.

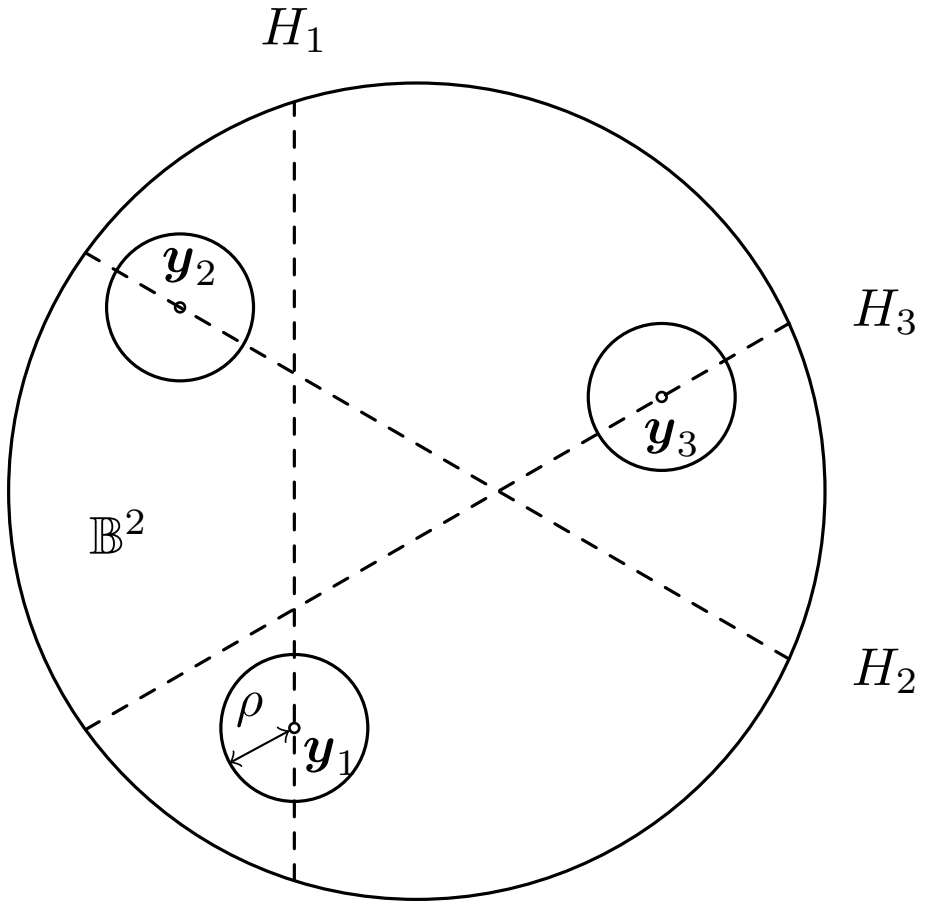
\begin{figure}[H]
\centering
\begin{tikzpicture}[x=2.7cm,y=2.7cm,font=\small,line cap=round]
 \draw[line width=.6pt] (0,0) circle[radius=1];
 \draw[dashed,line width=.6pt] (-.3,-.953939) -- (-.3,.953939);
 \draw[dashed,line width=.6pt] (-.811684,.584096) -- (.911684,-.410891);
 \draw[dashed,line width=.6pt] (-.811684,-.584096) -- (.911684,.410891);
 \draw[line width=.6pt] (-.3,-.58) circle[radius=.18];
 \draw[line width=.6pt] (-.58,.450333) circle[radius=.18];
 \draw[line width=.6pt] (.6,.230940) circle[radius=.18];
 \draw[line width=.5pt] (-.3,-.58) circle[radius=.012];
 \draw[line width=.5pt] (-.58,.450333) circle[radius=.012];
 \draw[line width=.5pt] (.6,.230940) circle[radius=.012];
 \draw[<->,line width=.4pt] (-.312,-.5865) -- (-.4584,-.6655);
 \node at (-.41,-.535) {\(\rho\)};
 \node at (-.21,-.65) {\(\vy_1\)};
 \node at (-.555,.550333) {\(\vy_2\)};
 \node at (.625,.130940) {\(\vy_3\)};
 \node[anchor=south] at (-.3,1.025) {\(H_1\)};
 \node[anchor=west] at (1.015,-.44) {\(H_2\)};
 \node[anchor=west] at (1.015,.44) {\(H_3\)};
 \node at (-.73,-.1) {\(\mathbb B^2\)};
\end{tikzpicture}
\caption{Schematic illustration of the parameter separation condition in $\mathbb{B}^2$.
Dashed lines represent the hyperplanes, and solid circles bound the disks.
On each disk \(B(\vy_j,\rho)\), every other atom agrees with a polynomial of degree at most \(s\), so the coefficient of the \(j\)-th atom can be estimated locally.}
\label{fig:parameter-separation}
\end{figure}

\leavevmode\phantomsection\label{not:ridge-hyperplane}For \(\vtheta=(\vxi,t)\in\Param\), let \(H_{\vxi,t}:=\{\vx\in\R^d:\vxi\cdot\vx=t\}\) be the associated hyperplane.
\leavevmode\phantomsection\label{not:indexed-hyperplane}For a finite parameter set \(\Lambda=\{\vtheta_1,\ldots,\vtheta_M\}\), abbreviate \(H_j=H_{\vxi_j,t_j}\).
We say that \(\Lambda\) satisfies the local hyperplane separation condition with radius \(\rho>0\) if, for every \(j\), there is a point \(\vy_j\in H_j\cap\Bd\) such that
\begin{equation}\label{eq:private-patch}
 B(\vy_j,\rho)\subset\Bd, \qquad B(\vy_j,\rho)\cap H_k=\varnothing\quad(k\ne j),
\end{equation}
and these balls are pairwise disjoint.
Figure~\ref{fig:parameter-separation} illustrates this condition for \(d=2\).
\Needspace{5cm}

\begin{theorem}
\label{thm:private-inverse}
Suppose that the scales satisfy \eqref{eq:geometric-scales}, and let \(\dict{j}\) index a finite network representation of each \(g_j\).
Suppose that \(\jointdict{j}\) satisfies the local hyperplane separation condition with radius \(\privrad{j}\), and that there are \(C_M,c_\rho>0\) and \(\nu,\kappa_0\ge0\) such that
\begin{equation*}
 \abs{\jointdict{j}}\le C_Mh_{j+1}^{-\nu}, \qquad \privrad{j}\ge c_\rho h_{j+1}^{\kappa_0}.
\end{equation*}
If
\begin{equation*}
 e_j\le Ah_j^\beta, \qquad \beta>\frac\nu2+\kappa_0\left(s+\frac d2\right),
\end{equation*}
then \(f\in\Vsp\).
\end{theorem}

The proof in Appendix~\ref{app:private-patches} gives \(\mingrameig{\jointdict{j}}\gtrsim_{d,s}\privrad{j}^{2s+d}\), so \eqref{eq:spectral-condition} holds with \(\eta=\kappa_0(s+d/2)\).
\section{Discussion and Conclusion}\label{sec:discussion-conclusion}

The direct estimate bounds the interpolation error in terms of the fill distance and the weighted variation norm.
The lower bound shows that boundary concentration determines the worst-case order under the boundary separation condition.
For the inverse problem, finite dual certificates control the growth of the interpolation norms.
The certificate argument permits singular Gram matrices and summable duality gaps.
Lower spectral bounds, obtained either directly or from local hyperplane separation, give explicit sufficient rates for these inverse conclusions.

Several questions remain open.
For the direct problem, the current matching lower bound relies on the boundary separation condition, while the optimal error for general node sets may also depend on finer geometric features, such as boundary distribution and local clustering.
A natural direction is to derive sharp direct bounds that account for the geometry of arbitrary node sets.
For the inverse problem, a closer study of the quantitative relation between interpolation data and the norms of dual certificates may lead to conditions that are easier to verify.
The geometry of atomic Banach spaces, characterizations of nonlinear approximation spaces, and the dual structure of interpolation problems also offer possible approaches to an inverse theory matching the direct estimates.

\appendix
\section*{Appendix}
\addcontentsline{toc}{section}{Appendix}
\renewcommand{\thesubsection}{\Alph{subsection}}
\renewcommand{\thesubsubsection}{\thesubsection.\arabic{subsubsection}}
\renewcommand{\theHsubsection}{appendix.\Alph{subsection}}
\counterwithin{equation}{subsection}
\counterwithin{definition}{subsection}
\counterwithin{theorem}{subsection}
\counterwithin{lemma}{subsection}
\counterwithin{proposition}{subsection}
\counterwithin{corollary}{subsection}
\counterwithin{example}{subsection}
\counterwithin{table}{subsection}

We now turn to the proofs and auxiliary results, which are organized into five appendices.
Appendix~\ref{app:weighted-network-auxiliary-results} establishes the properties and approximation estimates for weighted variation spaces.
Appendices~\ref{app:direct-proof} and~\ref{app:inverse-proofs} give the proofs of the direct and inverse results, respectively.
Appendix~\ref{app:private-patches} derives the Gram matrix spectral bound from local hyperplane separation, and Appendix~\ref{app:notation} collects the principal notation.

\Needspace{9\baselineskip}
\subsection{Weighted Variation Spaces}\label{app:weighted-network-auxiliary-results}
\leavevmode\phantomsection\label{not:unweighted-space}For comparison with the usual unweighted variation space, write
\begin{align*}
 \Vunw &:={\left\{ f=\int_{\Param}\rawatom{\vxi}{t}\dd{\nu}(\vtheta): \nu\in\Msp \right\}},\\
 \unwnorm{f} &:=\inf\left\{\tvnorm{\nu}: f=\int_{\Param}\rawatom{\vxi}{t}\dd{\nu}(\vtheta) \right\}.
\end{align*}
\subsubsection{Basic properties and comparison with the unweighted variation space}\label{app:basic-properties}
We first establish the compactness of the weighted dictionary, completeness of the space, and continuity of point evaluation.
We then compare the weighted and unweighted variation spaces.

\leavevmode\phantomsection\label{not:admissible-weight}More generally, a continuous nonnegative function \(w:\Param\to[0,\infty)\) is called an admissible weight if, after defining \(\rawatom{\vxi}{t}/w(\vtheta)\) to be zero whenever \(w(\vtheta)=0\), the resulting normalized atom depends continuously on \(\vtheta\) as an \(\Ltwo\)-valued function.
This notion and the corresponding construction of weighted dictionaries are given in \citep{devore2025weighted}.

\begin{lemma}\label{lem:cap-scale}
For every \(\vxi\in\Sph^{d-1}\) and \(-1\le t<1\), \leavevmode\phantomsection\label{not:cap-height}set \(\delta=1-t\).
Then
\begin{equation*}
 \norm[\Ltwo]{\rawatom{\vxi}{t}} \asymp_{d,s}\delta^{s+(d+1)/4}.
\end{equation*}
Consequently,
\begin{equation*}
 \norm[\Ltwo]{\atom{\vxi}{t}} \asymp_{d,s}\delta^{(2s+1)/4}.
\end{equation*}
In particular, the map \(\vtheta\mapsto\paramatom{\vtheta}\) maps the compact parameter space \(\Param\) continuously into \(\Ltwo\), and its image is a compact dictionary in \(\Ltwo\).
\end{lemma}
\begin{proof}
It suffices to prove the norm inequalities for \(\rawatom{\vxi}{t}\).
Since \(\vxi\in\Sph^{d-1}\), rotational invariance allows us to assume that \(\vxi=(0,\cdots,0,1)^\top\).
\leavevmode\phantomsection\label{not:cross-section-volume}Let \(\omega_{d-1}=\frac{\pi^{(d-1)/2}}{\Gamma((d+1)/2)}\) be the volume of the unit ball in \(\mathbb{R}^{d-1}\).
Then
\begin{align*}
 \left\lVert\rawatom{\vxi}{t}\right\rVert_{\Ltwo}^2 =&\int_{\Bd}(\vxi\cdot\vx-t)_+^{2s}\dd{\vx} =\int_{\Bd}(x_d-t)_+^{2s}\dd{\vx}\\
 =&\omega_{d-1}\int_t^1(x_d-t)^{2s}(1-x_d^2)^{(d-1)/2}\dd{x_d}\\
 =&\omega_{d-1}\delta^{2s+(d+1)/2} \int_0^1u^{2s}(1-u)^{(d-1)/2} (1+t+\delta u)^{(d-1)/2}\dd{u},
\end{align*}
where the last equality follows from the change of variables \(x_d=t+\delta u\).
Since \(1+t+\delta u=2-\delta(1-u)\le2\), we obtain \(\left\lVert\rawatom{\vxi}{t}\right\rVert_{\Ltwo}^2 \lesssim_{d,s}\delta^{2s+(d+1)/2}\).
For the lower bound, observe that
\begin{align*}
 \int_0^1u^{2s}(1-u)^{(d-1)/2}(1+t+\delta u)^{(d-1)/2}\dd{u} \ge&\int_{1/2}^{3/4}u^{2s}(1-u)^{(d-1)/2} (1+t+\delta u)^{(d-1)/2}\dd{u}\\
 \ge&\int_{1/2}^{3/4}2^{-2s-d+1}\dd{u} =2^{-(2s+d+1)}.
\end{align*}
Hence \(\left\lVert\rawatom{\vxi}{t}\right\rVert_{\Ltwo}^2 \gtrsim_{d,s}\delta^{2s+(d+1)/2}\), which proves the lemma.
\end{proof}

Lemma \ref{lem:cap-scale} implies that the synthesis operator \(\synthop:\Msp\to\Ltwo\) is bounded.
The integral defining \(\synth{\mu}\) is an \(\Ltwo\)-valued Bochner integral.
The measure representation and norm characterization in Section~\ref{sec:preliminaries} are standard properties of variation spaces generated by compact dictionaries \citep[Lemma~3]{siegelxu2023variation}; their counterpart for weighted ReLU dictionaries is given in \citep{devore2025weighted}.
Since \(\synthop\) is bounded, \(\ker\synthop\) is a closed subspace of \(\Msp\).
Hence \(\Vsp\) is isometrically isomorphic to the quotient space \(\Msp/\ker\synthop\), and is therefore a Banach space.
This also follows from the general completeness result for variation spaces \citep[Lemma~1]{siegelxu2023variation}.
Moreover,
\begin{equation*}
 \ltnorm{f} \le \sup_{\vtheta\in\Param}\norm[\Ltwo]{\paramatom{\vtheta}}\, \vnorm{f},
\end{equation*}
so \(\Vsp\hookrightarrow\Ltwo\) continuously.

Functions in \(\Vsp\) have continuous representatives in the interior.
Indeed, let \(K\Subset\Bd\) be arbitrary and write \leavevmode\phantomsection\label{not:compact-margin}\(\rho_K:=1-\sup_{\vx\in K}\abs{\vx}>0\).
If \(t>1-\rho_K\), then \(\vxi\cdot\vx-t<0\) for every \(\vx\in K\), and hence \(\atom{\vxi}{t}|_K=0\).
On the other hand, over the compact parameter set \(\Sph^{d-1}\times[-1,1-\rho_K]\), the denominator \((1-t)^{(d+2s)/4}\) is uniformly bounded away from zero, and \((\vxi,t,\vx)\mapsto\atom{\vxi}{t}(\vx)\) is continuous.
Therefore, the map
\begin{equation*}
 \Param\longrightarrow C(K),\qquad \vtheta\longmapsto\paramatom{\vtheta}|_K
\end{equation*}
is continuous and bounded.
The basic properties of the Bochner integral \citep[Chapter~II]{diesteluhl1977vector} give, \leavevmode\phantomsection\label{not:continuous-representative}for every \(\mu\in\Msp\),
\begin{equation*}
 F_\mu|_K:=\int_{\Param}\paramatom{\vtheta}|_K\dd{\mu}(\vtheta) \in C(K), \qquad \norm[C(K)]{F_\mu} \le \sup_{\vtheta\in\Param}\norm[C(K)]{\paramatom{\vtheta}}\, \tvnorm{\mu}.
\end{equation*}
This \(C(K)\)-valued integral agrees with \(\synth{\mu}|_K\) in \(L_2(K)\).
Since \(K\Subset\Bd\) is arbitrary, \(F_\mu\) gives a continuous representative of \(\synth{\mu}\) on \(\Bd\).
If \(\synth{\mu}=\synth{\nu}\) in \(\Ltwo\), then the continuous function \(F_\mu-F_\nu\) vanishes almost everywhere and therefore vanishes identically on \(\Bd\).
Thus, this continuous representative depends only on the \(\Ltwo\) equivalence class of the function.
Taking the infimum over all representing measures further gives
\begin{equation*}
 \norm[C(K)]{f} \le \sup_{\vtheta\in\Param}\norm[C(K)]{\paramatom{\vtheta}}\, \vnorm{f},
\end{equation*}
so \(\Vsp\hookrightarrow C(K)\) continuously for every \(K\Subset\Bd\).
In particular, evaluation at every point in the interior of the unit ball is a continuous linear functional on \(\Vsp\).

Proposition \ref{prop:unweighted-embedding} below gives a continuous embedding of the unweighted variation space into \(\Vsp\).
Example \ref{exm:strict-weighted-inclusion} then gives a one-dimensional example for which the inclusion is strict: a sum of ridge functions concentrating near the boundary has finite weighted variation norm and infinite unweighted variation norm.

\begin{proposition}
\label{prop:unweighted-embedding}
For every \(f\in\Vunw\),
\begin{equation*}
 \vnorm{f} \le 2^{(d+2s)/4}\unwnorm{f}.
\end{equation*}
Consequently, \(\Vunw\hookrightarrow\Vsp\).
\end{proposition}
\begin{proof}
Given an unweighted representing measure \(\nu\), define \(\dd{\mu}(\vtheta)=\wstar{t}\dd{\nu}(\vtheta)\).
Since \(0\le\wstar{t}\le 2^{(d+2s)/4}\), we have
\begin{equation*}
 \int_{\Param}\rawatom{\vxi}{t}\dd{\nu}(\vtheta) =\int_{\Param}\atom{\vxi}{t}\dd{\mu}(\vtheta), \qquad \tvnorm{\mu} \le 2^{(d+2s)/4}\tvnorm{\nu}.
\end{equation*}
Taking the infimum over all unweighted representations proves the claim.
\end{proof}

\begin{example}
\label{exm:strict-weighted-inclusion}
Let \(d=1\), so that \(\mathbb B^1=(-1,1)\), and set
\begin{equation*}
 t_j:=1-2^{-j},\qquad f(x):=\sum_{j=1}^{\infty}\frac1j(x-t_j)_+^s.
\end{equation*}
Then \(f\in\Vsp\), but \(f\notin\Vunw\).
\end{example}
\begin{proof}
When \(d=1\), \(\wstar{t_j}=2^{-j(2s+1)/4}\).
Therefore,
\begin{equation*}
 f=\sum_{j=1}^{\infty}c_j\atom{1}{t_j}, \qquad c_j:=\frac{2^{-j(2s+1)/4}}{j}, \qquad \sum_{j=1}^{\infty}\abs{c_j}<\infty,
\end{equation*}
and hence \(f\in\Vsp\).
On the other hand, in the sense of distributions,
\begin{equation*}
 \mathscr D^{s+1}f =s!\sum_{j=1}^{\infty}\frac1j\delta_{t_j}.
\end{equation*}
The total variation of the right-hand side on \((-1,1)\) is infinite.
If \(f\in\Vunw\), then applying \(\mathscr D^{s+1}\) to any representation by a measure of finite total variation would yield a finite signed Radon measure, contradicting the identity above.
\end{proof}

\subsubsection{Approximation of the weighted variation space by shallow networks}
We prove Theorem~\ref{thm:weighted-network-approximation} and Proposition~\ref{thm:weak-star-closure}.
For the unit ball of the unweighted space, an approximation bound \(O(n^{-1/2-(2s+1)/(2d)})\) in the uniform norm for shallow \(\relu^s\) networks was proved in \citep[Theorem~3]{siegel2025uniform}.
\leavevmode\phantomsection\label{not:weighted-parameter}For the weighted approximation estimate, let \(r=\sqrt{1-t}\), and parameterize the weighted dictionary by
\begin{equation*}
 \weightparam(\vxi,r)(\vx) :=
 \begin{cases}
 r^{-(d+2s)/2} \bigl(\vxi\cdot\vx-1+r^2\bigr)_+^s, &0<r\le\sqrt2,\\
 0,&r=0.
 \end{cases}
\end{equation*}
Then \(\weightparam(\vxi,r)=\atom{\vxi}{1-r^2}\).

\begin{lemma}
\label{lem:weighted-parameter-smoothness}
The map \(\weightparam:\Sph^{d-1}\times[0,\sqrt2]\longrightarrow\Ltwo\) belongs to \(\mathrm{Lip}_{\infty}^{s+1/2} \bigl(\Sph^{d-1}\times[0,\sqrt2];\Ltwo\bigr)\).
More precisely, in a fixed finite smooth atlas of \(\Sph^{d-1}\), for \(0\le\ell\le s\),
\begin{equation*}
 \norm[\Ltwo]{\mathscr D^{\ell}\weightparam(\vxi,r)} \lesssim_{d,s}r^{s+1/2-\ell}, \qquad 0<r\le1,
\end{equation*}
and the parameter derivatives of order \(s\) are \(1/2\)-H\"older continuous.
\end{lemma}

\begin{proof}
\leavevmode\phantomsection\label{not:cap-geometry}Set
\begin{equation*}
 q_{\vxi,r}(\vx):=\vxi\cdot\vx-1+r^2, \qquad \mathcal C_{\vxi,r}:= \{\vx\in\Bd:q_{\vxi,r}(\vx)>0\}.
\end{equation*}
Thus, \(\weightparam(\vxi,r)(\vx) =r^{-s-d/2}\bigl(q_{\vxi,r}(\vx)\bigr)_+^s\).
If \(0<r\le1\) and \(\vx\in\mathcal C_{\vxi,r}\), then \(0<q_{\vxi,r}(\vx)\le r^2\) and \(\abs{\vx-(\vxi\cdot\vx)\vxi}\le\sqrt2r\).
Moreover, integration over cross sections of the spherical cap gives \(\abs{\mathcal C_{\vxi,r}}\lesssim_dr^{d+1}\).

Fix a finite smooth atlas of \(\Sph^{d-1}\).
In a relatively compact parameter domain of one chart, \leavevmode\phantomsection\label{not:local-parameters}write \(\vxi=\vxi(\vect{u})\) and let \(\vect{z}=(\vect{u},r)\).
The first \(s+1\) derivatives of the finitely many coordinate maps have uniform bounds depending only on \(d,s\).
Since \(\partial_{u_j}\vxi(\vect{u})\cdot\vxi(\vect{u})=0\), for \(\vx\in\mathcal C_{\vxi,r}\),
\begin{align*}
 \abs{\partial_{u_j}q_{\vxi,r}(\vx)} &=\abs{\partial_{u_j}\vxi(\vect{u})\cdot \bigl(\vx-(\vxi\cdot\vx)\vxi\bigr)} \lesssim_d r,\\
 \abs{\partial_r q_{\vxi,r}(\vx)}&=2r.
\end{align*}
In addition, the absolute value of every parameter derivative of \(q_{\vxi,r}\) whose order lies between \(2\) and \(s+1\) is \(\lesssim_{d,s}1\).

Apply the product and chain rules repeatedly to \(r^{-s-d/2}(q_{\vxi,r})_+^s\).
For any parameter derivative of total order \(\ell\le s\), each term on \(\mathcal C_{\vxi,r}\) has the form
\begin{equation*}
 c\,r^{-s-d/2-a} (q_{\vxi,r})_+^{s-j} \prod_{\nu=1}^{j}\mathscr D_{\vect{z}}^{m_\nu}q_{\vxi,r}, \qquad a+\sum_{\nu=1}^{j}m_\nu=\ell, \qquad m_\nu\ge1.
\end{equation*}
Here, \(a\) is the number of derivatives falling on the normalization factor.
If \(n_1=\#\{\nu:m_\nu=1\}\), then \(0<q_{\vxi,r}\le r^2\) and the preceding estimates for the parameter derivatives show that the absolute value of this term is at most
\begin{align*}
 C_{d,s}r^{-s-d/2-a+2(s-j)+n_1} &=C_{d,s}r^{s-d/2-\ell+ \sum_{\{\nu:m_\nu\ge2\}}(m_\nu-2)}\le C_{d,s}r^{s-d/2-\ell}, \qquad 0<r\le1.
\end{align*}
It follows that
\begin{equation*}
 \abs{ \mathscr D_{\vect{z}}^{\ell} \left[r^{-s-d/2}(q_{\vxi,r})_+^s\right] } \lesssim_{d,s}r^{s-d/2-\ell} \boldsymbol{1}_{\mathcal C_{\vxi,r}}, \qquad 0\le\ell\le s.
\end{equation*}
Combining this with the volume estimate for the spherical cap yields
\begin{align*}
 \norm[\Ltwo]{ \mathscr D_{\vect{z}}^{\ell}\weightparam(\vxi,r) } &\lesssim_{d,s} r^{s-d/2-\ell} \abs{\mathcal C_{\vxi,r}}^{1/2}\lesssim_{d,s}r^{s+1/2-\ell}, \qquad 0\le\ell\le s.
\end{align*}
In particular, because \(s+1/2-\ell\ge1/2\), all these derivatives extend continuously by zero at \(r=0\).

We next verify the \(1/2\)-H\"older continuity of the derivatives of order \(s\).
Take \(\vect{z}=(\vect{u},r)\) and \(\vect{z}'=(\vect{u}',r')\) in the same coordinate chart, \leavevmode\phantomsection\label{not:parameter-distance}and set \(h:=\abs{\vect{z}-\vect{z}'}\).
Suppose first that \(h\le r/2\), so that \(r/2\le r'\le3r/2\).
Since \(\vxi(\vect{u})\) and \(\vxi(\vect{u}')\) are unit vectors, \(\bigl(\vxi(\vect{u}')-\vxi(\vect{u})\bigr) \cdot\vxi(\vect{u}) =-\frac12\abs{\vxi(\vect{u}')-\vxi(\vect{u})}^2\).
Together with the smoothness of the coordinate map, this identity gives \(\abs{q_{\vxi(\vect{u}),r}(\vx) -q_{\vxi(\vect{u}'),r'}(\vx)}\lesssim_d hr\) on the union of the two spherical caps.
Hence there is a constant \(C_d>0\), depending only on the coordinate chart, such that the symmetric difference of the two spherical caps is contained in
\begin{equation*}
 \left\{\vx\in\Bd: \abs{q_{\vxi(\vect{u}),r}(\vx)}\le C_dhr, \ \abs{\vx-(\vxi(\vect{u})\cdot\vx)\vxi(\vect{u})}\le C_dr \right\}.
\end{equation*}
Integration in the \(\vxi(\vect{u})\) direction gives
\begin{equation*}
 \abs{ \mathcal C_{\vxi(\vect{u}),r} \mathbin{\triangle} \mathcal C_{\vxi(\vect{u}'),r'} } \lesssim_dhr^d.
\end{equation*}
On the intersection of the two spherical caps, compare the terms obtained from the product rule above.
The normalization factors differ by at most \(C_{d,s}hr^{-s-d/2-a-1}\).
Also, the first derivatives of \(q\) differ by \(O_{d,s}(h)\), and its derivatives of order at least two differ by \(O_{d,s}(h)\).
Finally, for \(k\ge1\),
\begin{equation*}
 \abs{q_{\vxi,r}(\vx)^k-q_{\vxi',r'}(\vx)^k} \lesssim_{d,s}hr^{2k-1}
\end{equation*}
on this intersection, since both positive values are \(O(r^2)\) and their difference is \(O(hr)\).
Replacing one factor at a time in a term of total derivative order \(s\) therefore changes its value by at most \(C_{d,s}hr^{-d/2-1}\).
This gives the pointwise estimate
\begin{equation*}
 \abs{ \mathscr D_{\vect{z}}^s\weightparam(\vect{z})(\vx) -\mathscr D_{\vect{z}}^s\weightparam(\vect{z}')(\vx) } \lesssim_{d,s}hr^{-d/2-1}.
\end{equation*}
The volume estimate for spherical caps and \(r/2\le r'\le3r/2\) imply \(\abs{\mathcal C_{\vxi(\vect{u}),r} \cup\mathcal C_{\vxi(\vect{u}'),r'}}\lesssim_dr^{d+1}\).
Therefore,
\begin{equation*}
 \norm[\Ltwo]{ \boldsymbol{1}_{ \mathcal C_{\vxi(\vect{u}),r} \cap\mathcal C_{\vxi(\vect{u}'),r'}} \left( \mathscr D_{\vect{z}}^s\weightparam(\vect{z}) -\mathscr D_{\vect{z}}^s\weightparam(\vect{z}') \right) } \lesssim_{d,s}hr^{-1/2} \lesssim_{d,s}h^{1/2}.
\end{equation*}
Each term on the symmetric difference has magnitude \(\lesssim_{d,s}r^{-d/2}\), and hence its \(L_2\) norm is at most a constant multiple of \(r^{-d/2}(hr^d)^{1/2}=h^{1/2}\).
If \(h>r/2\), then \(r'\le r+h\le3h\), and the estimate already proved for derivatives of order \(s\) gives
\begin{align*}
 \norm[\Ltwo]{ \mathscr D_{\vect{z}}^s\weightparam(\vect{z}) -\mathscr D_{\vect{z}}^s\weightparam(\vect{z}') } &\lesssim_{d,s}r^{1/2}+(r')^{1/2}\lesssim_{d,s}h^{1/2}.
\end{align*}
This proves that, when \(0\le r,r'\le1\), the parameter derivatives of order \(s\) satisfy
\begin{equation*}
 \norm[\Ltwo]{ \mathscr D_{\vect{z}}^s\weightparam(\vect{z}) -\mathscr D_{\vect{z}}^s\weightparam(\vect{z}') } \lesssim_{d,s}\abs{\vect{z}-\vect{z}'}^{1/2}.
\end{equation*}

On the parameter region \(1/2\le r\le\sqrt2\), the normalization factor \(r^{-s-d/2}\) and its first \(s+1\) derivatives are uniformly bounded.
The conclusion then reduces to the parameter smoothness of the usual \(\relu^s\) dictionary \citep[Lemma~1, $p=2$]{siegelxu2024sharp}.
Combining the estimates above over the finite atlas of \(\Sph^{d-1}\) yields
\begin{equation*}
 \weightparam\in \mathrm{Lip}_{\infty}^{s+1/2} \bigl(\Sph^{d-1}\times[0,\sqrt2];\Ltwo\bigr).
\end{equation*}
\end{proof}

\begin{proof}[Proof of Theorem \ref{thm:weighted-network-approximation}]
Lemma \ref{lem:weighted-parameter-smoothness} shows that the normalized weighted dictionary is parameterized by a compact \(d\)-dimensional manifold with boundary, and the parameterization has smoothness of order \(s+1/2\).
Applying the approximation theorem for smoothly parameterized dictionaries \citep[Theorem~2]{siegelxu2024sharp} in the Hilbert space \(\Ltwo\) gives
\begin{equation*}
 \bestnerr{f} \lesssim_{d,s}n^{-1/2-(s+1/2)/d}\vnorm{f}.
\end{equation*}
That theorem also provides \(\ell_1\) control of the normalized atomic coefficients.
Rewriting the normalized atoms in terms of the original \(\relu^s\) atoms gives the weighted atomic cost estimate.
\end{proof}

\begin{proof}[Proof of Proposition \ref{thm:weak-star-closure}]
Let \(L=\liminf_{n\to\infty}\vnorm{g_n}\), and choose a subsequence, still denoted by \((g_n)\), such that \(\vnorm{g_n}\to L\).
By the definition of the norm, one can choose \(\mu_n\in\Msp\) satisfying
\begin{equation*}
 \synth{\mu_n}=g_n, \qquad \tvnorm{\mu_n}\le\vnorm{g_n}+n^{-1}.
\end{equation*}
The norm definition also gives \(\vnorm{g_n}\le\tvnorm{\mu_n}\), and hence \(\tvnorm{\mu_n}\to L\).
Thus, \((\mu_n)\) is uniformly bounded in \(C(\Param)^*=\Msp\), where the duality is given by the Riesz representation theorem \citep[Corollary~7.18]{folland1999real}.
Since \(\Param\) is a compact metrizable space, \(C(\Param)\) is separable, and the Banach--Alaoglu theorem \citep[Theorem~5.18]{folland1999real} yields a subsequence and a measure \(\mu\in\Msp\) such that \(\mu_n\overset{*}{\rightharpoonup}\mu\).

\leavevmode\phantomsection\label{not:dual-test}Now take any \(v\in\Ltwo\).
By Lemma \ref{lem:cap-scale}, the function
\begin{equation*}
 H_v(\vtheta):=\ip{\paramatom{\vtheta}}{v}_{\Ltwo}
\end{equation*}
belongs to \(C(\Param)\).
Therefore,
\begin{align*}
 \ip{g_n}{v}_{\Ltwo} &=\int_{\Param}H_v(\vtheta)\dd{\mu_n}(\vtheta)\longrightarrow \int_{\Param}H_v(\vtheta)\dd{\mu}(\vtheta) =\ip{\synth{\mu}}{v}_{\Ltwo}.
\end{align*}
This shows that \(g_n\rightharpoonup\synth{\mu}\) in \(\Ltwo\).
On the other hand, the strong convergence in the hypothesis implies \(g_n\rightharpoonup f\).
Uniqueness of the weak limit gives \(f=\synth{\mu}\), and hence \(f\in\Vsp\).
Finally, the weak-star lower semicontinuity of the total variation norm yields
\begin{equation*}
 \vnorm{f} \le\tvnorm{\mu} \le\liminf_{n\to\infty}\tvnorm{\mu_n} =L.
\end{equation*}
This completes the proof.
\end{proof}

\begin{corollary}
\label{cor:norm-measure-consistency}
Let \(f\in\Vsp\), let the node family \(X_j\) satisfy the sampling conditions of Theorem \ref{thm:direct}, and suppose that \(h_j\to0\).
For any choice \(g_j=\ixjf\),
\begin{equation*}
 g_j\to f\quad\text{in }\Ltwo, \qquad \lim_{j\to\infty}\vnorm{g_j}=\vnorm{f}.
\end{equation*}
Furthermore, if \(\eps_j\downarrow0\) and \(\mu_j\in\Msp\) is chosen so that
\begin{equation*}
 \synth{\mu_j}=g_j, \qquad \tvnorm{\mu_j}\le\vnorm{g_j}+\eps_j,
\end{equation*}
then every weak-star cluster point \(\mu\) of \((\mu_j)\) satisfies
\begin{equation*}
 \synth{\mu}=f, \qquad \tvnorm{\mu}=\vnorm{f}.
\end{equation*}
If \(f\) has a unique representing measure of minimum total variation, then the entire sequence satisfies \(\mu_j\overset{*}{\rightharpoonup}\mu\).
\end{corollary}

\begin{proof}
Theorem \ref{thm:direct} first gives \(g_j\to f\) in \(\Ltwo\).
Since \(f\) itself satisfies all the interpolation conditions, minimality gives
\begin{equation*}
 \vnorm{g_j}\le\vnorm{f}.
\end{equation*}
On the other hand, the lower semicontinuity in Proposition \ref{thm:weak-star-closure} shows that
\begin{equation*}
 \vnorm{f} \le\liminf_{j\to\infty}\vnorm{g_j}.
\end{equation*}
Combining these two inequalities gives \(\vnorm{g_j}\to\vnorm{f}\).

It follows from this convergence and \(\eps_j\downarrow0\) that the total variations of the measures \((\mu_j)\) are uniformly bounded.
Take any weak-star convergent subsequence and, after relabeling, write \(\mu_j\overset{*}{\rightharpoonup}\mu\).
As in the proof of Proposition \ref{thm:weak-star-closure}, for arbitrary \(v\in\Ltwo\), testing against \(H_v(\vtheta)=\ip{\paramatom{\vtheta}}{v}_{\Ltwo}\) gives \(\synth{\mu_j}\rightharpoonup\synth{\mu}\) in \(\Ltwo\).
Since \(\synth{\mu_j}=g_j\to f\) in \(\Ltwo\), we have \(\synth{\mu}=f\).
Finally, weak-star lower semicontinuity of the total variation norm yields
\begin{equation*}
 \vnorm{f} \le\tvnorm{\mu} \le\liminf_{j\to\infty}\tvnorm{\mu_j} \le\vnorm{f},
\end{equation*}
so \(\tvnorm{\mu}=\vnorm{f}\).
If the representing measure of minimum total variation is unique, then all weak-star convergent subsequences of \((\mu_j)\) have the same limit.
Since the weak-star topology is metrizable on balls of bounded total variation, the entire sequence converges weak-star to that limit.
This proves the corollary.
\end{proof}

\subsection{Proofs of Theorems \ref{thm:direct} and \ref{thm:direct-lower}}\label{app:direct-proof}
We construct a locally stable reconstruction operator that depends only on the nodal values and reproduces polynomials of degree at most \(s\).
\leavevmode\phantomsection\label{not:polynomial-space}Let \(\Polys=\Pi_s(\R^d)\) denote the space of polynomials on \(\R^d\) of total degree at most \(s\).
\begin{proposition}
\label{prop:local-reconstruction}
There exist constants \(h_0,R_s,C_s,N_{d,s,\kappa}>0\), depending only on \(d,s,\kappa\), with the following property.
Whenever \(0<h<h_0\) and \eqref{eq:fill-distance} and \eqref{eq:interior-margin} hold, \(\Bd\) admits a partition into a finite family \(\mathscr T_h\) of measurable cells of diameter at most \(C_sh\).
Using only values at the nodes in \(X_h\), one can construct a piecewise polynomial operator \(\sampling{h}\) such that the partition and the operator satisfy:
\begin{enumerate}
 \item \leavevmode\phantomsection\label{not:local-nodes}For every \(K\in\mathscr T_h\), there is a local node set \(Y_K\subset \{\vy\in X_h:\dist(\vy,K)\le R_sh\}\) such that \(\#Y_K=\dim\Polys\).

 \item The operator \(\sampling{h}\) reproduces every polynomial in \(\Polys\), that is,
 \[
 \sampling{h}p=p\qquad(p\in\Polys).
 \]
 In particular, if \(v|_{X_h}=0\), then \(\sampling{h}v=0\) on \(\Bd\).

 \item For every \(K\in\mathscr T_h\) and all data \(v\) on \(X_h\), the following local stability estimate holds:
 \[
 \norm[L_\infty(K)]{\sampling{h}v} \le C_s\max_{\vy\in Y_K}\abs{v(\vy)}.
 \]

 \item \leavevmode\phantomsection\label{not:enlarged-cell}Let
 \[
 K_{R_s}^{+}:=\{\vx\in\R^d:\dist(\vx,K)\le R_sh\}.
 \]
 Then the family of enlarged cells \(\{K_{R_s}^{+}:K\in\mathscr T_h\}\) has uniformly bounded overlap.
 More precisely,
 \[
 \sup_{\vx\in\R^d} \#\{K\in\mathscr T_h:\vx\in K_{R_s}^{+}\} =\sup_{\vx\in\R^d} \sum_{K\in\mathscr T_h}\mathbf 1_{K_{R_s}^{+}}(\vx) \le N_{d,s,\kappa}.
 \]
\end{enumerate}
\end{proposition}
\begin{proof}
We first construct \(\mathscr{T}_h\) and the local node sets \(Y_K\).
Consider half-open cubes of side length \(h\) and the family of their intersections with \(\Bd\):
\begin{align*}
 \mathscr{T}_h=\left\{h(\vect{k}+[0,1)^d)\cap\Bd:\vect{k}\in\mathbb{Z}^d, \ h(\vect{k}+[0,1)^d)\cap\Bd\ne\varnothing\right\}.
\end{align*}
Since each \(K\in\mathscr T_h\) is contained in a cube of side length \(h\), \(\operatorname{diam}K\le\sqrt{d}h\).
\leavevmode\phantomsection\label{not:unisolvent-template}Let \(M=\dim\Pi_s\), and choose a monomial basis \(q_1,\cdots,q_M\) for \(\Pi_s\).
Choose an arbitrary unisolvent point set \(\widehat{Y}=\left\{\widehat{\vy}_1,\cdots,\widehat{\vy}_M\right\}\) in \(\Bd\).
\leavevmode\phantomsection\label{not:vandermonde}For any labeled point configuration \(U=\{\vu_1,\ldots,\vu_M\}\), write
\begin{align*}
 \matr{V}(U):=\bigl(q_\ell(\vu_j)\bigr)_{j,\ell=1}^{M}.
\end{align*}
The unisolvency of \(\widehat{Y}\) is equivalent to the invertibility of the Vandermonde matrix \(\matr{V}(\widehat{Y})\).
Since the Vandermonde determinant depends continuously on the nodes, \leavevmode\phantomsection\label{not:unisolvency-constants}there exists \(\varepsilon_s>0\) such that whenever \(|\vu_j-\widehat{\vy}_j|\le\varepsilon_s\), the new point configuration \(U=\{\vu_1,\cdots,\vu_M\}\) remains unisolvent and
\begin{align*}
 B_s:=\sup\left\{\lVert\matr{V}(U)^{-1}\rVert_{\ell_\infty\to\ell_\infty}: |\vu_j-\widehat{\vy}_j|\le\varepsilon_s,\ 1\le j\le M\right\}<\infty.
\end{align*}
\leavevmode\phantomsection\label{not:reconstruction-scales}Choose a sufficiently large constant \(L_s\) such that
\begin{align*}
 L_s>\frac{2}{\varepsilon_s},\qquad L_s\cdot\min_{i\neq j}|\widehat{\vy}_i-\widehat{\vy}_j|>4,
\end{align*}
and define
\begin{align*}
 A_s:=\kappa+L_s+2,\qquad D_s:=A_s+1,\qquad h_0=\frac{1}{2D_s}.
\end{align*}
\leavevmode\phantomsection\label{not:cell-center}For each \(K\in\mathscr{T}_h\), choose \(\vx_K\in K\) and set
\begin{align*}
 \vc_K=\begin{cases} \vx_K,&\ 1-|\vx_K|\ge A_sh\\
 (1-D_sh)\frac{\vx_K}{|\vx_K|},&\ \text{otherwise}
 \end{cases}.
\end{align*}
If \(\vc_K=\vx_K\), then \(\operatorname{dist}(\vc_K,K)=0\).
Otherwise, since \(|\vx_K|<1\),
\begin{align*}
 \operatorname{dist}(\vc_K,K) \le|\vc_K-\vx_K| =|\vx_K|-(1-D_sh)<D_sh.
\end{align*}
In the first case, \(|\vc_K|\le1-A_sh\); in the second, \(|\vc_K|=1-D_sh\).
\leavevmode\phantomsection\label{not:template-nodes}Therefore, for \(\vz_{K,j}=\vc_K+L_sh\widehat{\vy}_j\),
\begin{align*}
 |\vz_{K,j}| \le|\vc_K|+L_sh<1,
\end{align*}
so \(\vz_{K,j}\in\Bd\).
By the definition of the fill distance, \leavevmode\phantomsection\label{not:selected-nodes}there exists \(\vy_{K,j}\in X_h\) such that \(|\vy_{K,j}-\vz_{K,j}|\le h\).
If \(\vy_{K,i}=\vy_{K,j}\), then
\begin{align*}
 |\vz_{K,i}-\vz_{K,j}| \le|\vz_{K,i}-\vy_{K,i}|+|\vy_{K,j}-\vz_{K,j}|\le2h,
\end{align*}
which contradicts \(|\vz_{K,i}-\vz_{K,j}|=L_s|\widehat{\vy}_i-\widehat{\vy}_j|h>4h\).
Thus the chosen nodes are distinct, and \(Y_K=\{\vy_{K,1},\cdots,\vy_{K,M}\}\) is the required local node set.
Indeed,
\begin{align*}
 \operatorname{dist}(\vy_{K,j},K) \le&|\vy_{K,j}-\vz_{K,j}|+|\vz_{K,j}-\vc_K|+\operatorname{dist}(\vc_K,K)\\
 \le&h+L_sh+D_sh=(1+L_s+D_s)h.
\end{align*}
Thus, setting \(R_s=(1+L_s+D_s)\) gives \(Y_K\subset\{\vy\in X_h:\operatorname{dist}(\vy,K)\le R_sh\}\).

We next construct the piecewise polynomial operator \(\sampling{h}\).
\leavevmode\phantomsection\label{not:scaled-nodes}Let \(\vu_{K,j}=\frac{\vy_{K,j}-\vc_K}{L_sh}\).
Then \(|\vu_{K,j}-\widehat{\vy}_j|\le\frac{1}{L_s}<\varepsilon_s\), so \(U_K=\{\vu_{K,1},\cdots,\vu_{K,M}\}\) is unisolvent.
Since affine transformations preserve unisolvency, \(Y_K\) is also unisolvent.
For arbitrary data \(v\) on \(X_h\), let \(p_K\in\Polys\) be the unique polynomial that interpolates \(v\) on \(Y_K\).
Define \(\sampling{h}v\) by setting \(\sampling{h}v|_K=p_K|_K\) for each \(K\in\mathscr T_h\).
We verify that the resulting operator \(\sampling{h}\) has the stated properties.
For any \(p\in\Pi_s\), uniqueness implies that \(\sampling{h}p\) agrees with \(p\) on each \(K\), and hence \(\sampling{h}p=p\) throughout \(\Bd\).
If \(v|_{X_h}=0\), then the data for every local interpolation problem vanish, so uniqueness gives \(\sampling{h}v=0\).
\leavevmode\phantomsection\label{not:scaled-polynomial}To prove property 3, define the scaled interpolation polynomial
\begin{align*}
 \widehat{\mathscr P}_Kv(\vu):=p_K(\vc_K+L_sh\vu) =\sum_{\ell=1}^{M}a_{K,\ell}q_\ell(\vu),
\end{align*}
\leavevmode\phantomsection\label{not:local-vectors}and let its coefficient vector and nodal data vector be, respectively,
\begin{align*}
 \vect{a}_K&=(a_{K,1},\ldots,a_{K,M})^{\mathsf T},\\
 \vect{b}_K&=(v(\vy_{K,1}),\ldots,v(\vy_{K,M}))^{\mathsf T}.
\end{align*}
\leavevmode\phantomsection\label{not:local-vandermonde}Set
\begin{align*}
 \matr{V}_K:=\matr{V}(U_K)=\bigl(q_\ell(\vu_{K,j})\bigr)_{j,\ell=1}^{M}.
\end{align*}
The interpolation conditions \(\widehat{\mathscr P}_Kv(\vu_{K,j})=v(\vy_{K,j})\) give
\begin{align*}
 \matr{V}_K\vect{a}_K=\vect{b}_K.
\end{align*}
By the bound on \(\matr V(U_K)^{-1}\) in the definition of \(B_s\),
\begin{align*}
 \lVert\vect{a}_K\rVert_{\ell_\infty} \le B_s\lVert\vect{b}_K\rVert_{\ell_\infty} =B_s\max_{\vy\in Y_K}|v(\vy)|.
\end{align*}
\leavevmode\phantomsection\label{not:scaled-cell-radius}On the other hand, for every \(\vx\in K\), using \(\vx_K\in K\), \(\operatorname{diam}K\le\sqrt d h\), and \(|\vx_K-\vc_K|\le D_sh\), we have
\begin{align*}
 \left|\frac{\vx-\vc_K}{L_sh}\right| &\le\frac{|\vx-\vx_K|+|\vx_K-\vc_K|}{L_sh} \le\frac{\sqrt d+D_s}{L_s}=: \beta_{d,s,\kappa}.
\end{align*}
\leavevmode\phantomsection\label{not:basis-bound}Because \(q_1,\ldots,q_M\) are fixed polynomials, \(C_q:=\sup_{|\vu|\le\beta_{d,s,\kappa}} \sum_{\ell=1}^{M}|q_\ell(\vu)|\) is finite.
Letting \(\vu=(\vx-\vc_K)/(L_sh)\), for every \(\vx\in K\),
\begin{align*}
 |(\sampling{h}v)(\vx)| &=|\widehat{\mathscr P}_Kv(\vu)| \le\sum_{\ell=1}^{M}|a_{K,\ell}|\,|q_\ell(\vu)|\\
 &\le C_q\lVert\vect{a}_K\rVert_{\ell_\infty} \le C_qB_s\max_{\vy\in Y_K}|v(\vy)|.
\end{align*}
Thus, taking \(C_s=C_qB_s\) gives
\begin{align*}
 \norm[L_\infty(K)]{\sampling{h}v} \le C_s\max_{\vy\in Y_K}|v(\vy)|.
\end{align*}

It remains to verify the bounded overlap asserted in property 4.
Fix \(\vx\in\R^d\) and suppose that \(\vx\in K_{R_s}^+\).
\leavevmode\phantomsection\label{not:parent-cube}Let \(\widetilde{K}=h(\bm{k}+[0,1)^d)\) be the parent cube associated with this enlarged cell.
We claim that \(\widetilde{K}\subset B(\vx,(R_s+\sqrt{d}+1)h)\).
Indeed, for every \(\vy\in\widetilde{K}\),
\begin{align*}
 |\vx-\vy|\le|\vx-\vy_K|+|\vy_K-\vy| \le (R_s+1)h+\sqrt{d}h=(R_s+\sqrt{d}+1)h,
\end{align*}
where \(\vy_K\) is a point of \(K\) satisfying \(|\vx-\vy_K|\le(R_s+1)h\).
Since this inclusion holds for every enlarged cell containing \(\vx\), while the parent cubes are disjoint, comparison of volumes yields
\begin{align*}
 \#\{K\in\mathscr{T}_h:\vx\in K_{R_s}^+\}h^d \le|B(\vx,(R_s+\sqrt{d}+1)h)| =\frac{\pi^{d/2}}{\Gamma(d/2+1)}(R_s+\sqrt{d}+1)^dh^d.
\end{align*}
Canceling \(h^d\) and setting \(N_{d,s,\kappa}=\frac{\pi^{d/2}}{\Gamma(d/2+1)}(R_s+\sqrt{d}+1)^d\) proves the proposition.
\end{proof}

By the polynomial reproduction property in Proposition \ref{prop:local-reconstruction}, reconstruction is exact on a cell whenever the ridge atom agrees with a polynomial of degree at most \(s\) throughout the neighborhood used for reconstruction.
For a ridge atom with associated hyperplane \(H_{\vxi,t}\), the reconstruction error can occur only when \(K_{R_s}^{+}\cap H_{\vxi,t}\ne\varnothing\) for a cell \(K\); for such cells, \(\operatorname{dist}(K,H_{\vxi,t})\le R_sh\).
We first establish a uniform estimate for a single normalized atom and then extend it to the entire weighted variation space through the measure representation.

\begin{lemma}\label{lem:operator-error}
Under the assumptions of Proposition \ref{prop:local-reconstruction}, for every \(u\in\Vsp\),
\begin{align*}
 \left\lVert u-\sampling{h}u\right\rVert_{L_2(\Bd)} \lesssim_{d,s,\kappa}h^{(2s+1)/4}\vnorm{u}.
\end{align*}
In particular, if \(u|_{X_h}=0\), then
\begin{align*}
 \lVert u\rVert_{L_2(\Bd)}\lesssim_{d,s,\kappa}h^{(2s+1)/4}\vnorm{u}.
\end{align*}
\end{lemma}
\begin{proof}
We first estimate the error for a single atom.
Fix \(\vtheta=(\vxi,t)\in\Param\), let \(\delta=1-t\), and let \(\paramatom{\vtheta}\) be the normalized atom defined in \eqref{eq:normalized-cap-atom}.
We prove the following estimate for the reconstruction error of \(\sampling{h}\) on \(\paramatom{\vtheta}\):
\begin{align}\label{eq:single-atom-interpolation-error}
 \left\lVert\paramatom{\vtheta}-\sampling{h}\paramatom{\vtheta}\right\rVert_{L_2(\Bd)} \lesssim_{d,s,\kappa}h^{(2s+1)/4}.
\end{align}

\paragraph{Case 1:} \(\delta<\kappa h\).
In this case, \(\vxi\cdot\vx_i-t\le|\vx_i|-(1-\delta) \le1-\kappa h-1+\delta=\delta-\kappa h<0\), so \(\paramatom{\vtheta}(\vx_i)=0\) at every node, and consequently \(\sampling{h}\paramatom{\vtheta}=0\).
Lemma \ref{lem:cap-scale} now gives \eqref{eq:single-atom-interpolation-error}.

\paragraph{Case 2:} \(\delta\ge\kappa h\).
Consider the hyperplane \(H_{\vxi,t}\) defined in Section~\ref{sec:inverse}.
If \(K_{R_s}^+\cap H_{\vxi,t}=\varnothing\), then the ridge atom \(\paramatom{\vtheta}\) is a polynomial on \(K_{R_s}^+\), so \(\sampling{h}\paramatom{\vtheta}|_K=\paramatom{\vtheta}|_K\), and hence \(\lVert\paramatom{\vtheta}-\sampling{h}\paramatom{\vtheta}\rVert_{L_2(K)}=0\).
Thus the error occurs only on cells \(K\) for which \(K_{R_s}^+\cap H_{\vxi,t}\ne\varnothing\).

Fix such a cell \(K\), and choose \(\vz\in K_{R_s}^+\cap H_{\vxi,t}\).
For every \(\vx\in K\) and \(\vy\in Y_K\), the diameter of \(K\) and the locality of \(Y_K\) imply \(|\vx-\vz|+|\vy-\vz|\lesssim_{d,s,\kappa}h\), and therefore
\begin{align*}
 |\paramatom{\vtheta}(\vx)| +|\paramatom{\vtheta}(\vy)| \lesssim_{d,s,\kappa}h^s\delta^{-(d+2s)/4}.
\end{align*}
The local stability of \(\sampling{h}\) gives
\begin{align*}
 \left\lVert\paramatom{\vtheta}-\sampling{h}\paramatom{\vtheta}\right\rVert_{L_2(K)} &\le|K|^{1/2}\left\lVert\paramatom{\vtheta}-\sampling{h}\paramatom{\vtheta}\right\rVert_{L_\infty(K)}\lesssim_{d,s,\kappa}h^{s+d/2}\delta^{-(d+2s)/4}.
\end{align*}
Since the interiors of the cells are disjoint,
\begin{align}\label{eq:single-atom-interpolation-error-2}
 \left\lVert\paramatom{\vtheta}-\sampling{h}\paramatom{\vtheta}\right\rVert_{L_2(\Bd)}^2 &=\sum_{\substack{K\in\mathscr T_h:\\
 K_{R_s}^+\cap H_{\vxi,t}\ne\varnothing}} \left\lVert\paramatom{\vtheta}-\sampling{h}\paramatom{\vtheta}\right\rVert_{L_2(K)}^2\lesssim_{d,s,\kappa}N_h(\vtheta)h^{2s+d}\delta^{-(d+2s)/2},
\end{align}
\leavevmode\phantomsection\label{not:affected-cells}where \(N_h(\vtheta)\) denotes the number of these cells.
We estimate \(N_h(\vtheta)\) by comparing volumes.
\leavevmode\phantomsection\label{not:slab}Define the region within distance \(Ah\) of the hyperplane by
\begin{align*}
 U_{h,A,C}(\vtheta):=\left\{\vx\in\mathbb{B}_{1+Ch}^d:|\vxi\cdot\vx-t|\le Ah\right\},
\end{align*}
where \(A,C\) are constants to be determined.
We choose \(A,C\) so that the cubes containing these cells lie in this slab.
More precisely, if \(K_{R_s}^+\cap H_{\vxi,t}\neq\varnothing\) and \(\vx\in\widetilde{K}=h(\bm{k}+[0,1)^d)\), then \(|\vxi\cdot\vx-t|\le(R_s+\sqrt{d}+1)h\).
Also, \(K\subset\Bd\) gives \(|\vx|\le1+(\sqrt{d}+1)h\).
Thus, taking \(A=R_s+\sqrt{d}+1\) and \(C=\sqrt{d}+1\) ensures that \(\widetilde{K}\subset U_{h,A,C}(\vtheta)\).
Let \(R=1+Ch\).
Slicing in the \(\vxi\) direction gives
\begin{align*}
 |U_{h,A,C}(\vtheta)| &=\omega_{d-1}\int_{[t-Ah,t+Ah]\cap[-R,R]} (R^2-z^2)^{(d-1)/2}\dd{z}.
\end{align*}
When \(0<h\le1\) and \(|z-t|\le Ah\), \(R^2-z^2\lesssim_{A,C}1-t^2+h\).
The length of the interval of integration is at most \(2Ah\), and hence
\begin{align*}
 |U_{h,A,C}(\vtheta)|\lesssim_{d,A,C}h(1-t^2+h)^{(d-1)/2}.
\end{align*}
The disjointness of the parent cubes gives
\begin{align*}
 N_h(\vtheta)&\le|U_{h,A,C}(\vtheta)|h^{-d}\\
 &\lesssim_{d,s,\kappa}(1-t^2+h)^{(d-1)/2}h^{1-d}\\
 &\lesssim_{d,s,\kappa}h^{1-d}(1-t^2)^{(d-1)/2}+h^{-(d-1)/2}.
\end{align*}
Substituting this into \eqref{eq:single-atom-interpolation-error-2} and using \(\delta\ge\kappa h\), we obtain
\begin{align*}
 \left\lVert\paramatom{\vtheta}-\sampling{h}\paramatom{\vtheta}\right\rVert_{L_2(\Bd)}^2 &\lesssim_{d,s,\kappa}h^{2s+1}\delta^{-(d+2s)/2}(1-t^2)^{(d-1)/2} +h^{2s+(d+1)/2}\delta^{-(d+2s)/2}\\
 &\lesssim_{d,s,\kappa}h^{2s+1}\delta^{-(2s+1)/2} +h^{2s+(d+1)/2}\delta^{-(d+2s)/2}\\
 &\lesssim_{d,s,\kappa}h^{(2s+1)/2}.
\end{align*}
Taking square roots proves \eqref{eq:single-atom-interpolation-error}.

Now take any \(u\in\Vsp\) and a representation \(u=\synth\mu\).
\leavevmode\phantomsection\label{not:lagrange-basis}For each cell \(K\), let \(\{\ell_{K,\vy}\}_{\vy\in Y_K}\) denote the Lagrange basis associated with the nodes \(Y_K\).
Using the continuous representative constructed in Appendix~\ref{app:basic-properties}, we have \(u(\vy)=\int_{\Param}\paramatom{\vtheta}(\vy)\dd{\mu}(\vtheta)\) for each \(\vy\in Y_K\).
Therefore,
\begin{align*}
 \sampling{h}u|_K &=\sum_{\vy\in Y_K}u(\vy)\ell_{K,\vy} =\int_{\Param}\sampling{h}\paramatom{\vtheta}|_K\dd{\mu}(\vtheta).
\end{align*}
Since \(\mathscr T_h\) is finite, this identity requires only interchanging finite sums with scalar integrals, and hence
\begin{align*}
 \sampling{h}u &=\sum_{K\in\mathscr T_h}\mathbf 1_K\sampling{h}u|_K\\
 &=\int_{\Param}\sum_{K\in\mathscr T_h}\mathbf 1_K \sampling{h}\paramatom{\vtheta}|_K\dd{\mu}(\vtheta)\\
 &=\int_{\Param}\sampling{h}\paramatom{\vtheta}\dd{\mu}(\vtheta).
\end{align*}
For fixed \(h\), the reconstruction \(\sampling{h}\paramatom{\vtheta}\) is a linear combination of finitely many fixed \(L_2\) functions with coefficients that depend continuously on the parameter, so the Bochner integral above is well defined.
\leavevmode\phantomsection\label{not:polar-density}Write \(\dd\mu=\zeta\,\dd|\mu|\), where \(|\zeta|=1\) \(|\mu|\)-almost everywhere.
Using \eqref{eq:single-atom-interpolation-error} and the triangle inequality for the Bochner integral \citep{diesteluhl1977vector}, we obtain
\begin{align*}
 \lVert u-\sampling{h}u\rVert_{L_2(\Bd)} &=\left\lVert\int_{\Param}\left(\paramatom{\vtheta}-\sampling{h}\paramatom{\vtheta}\right)\dd\mu(\vtheta)\right\rVert_{L_2(\Bd)}\\
 &\le\int_{\Param}\left\lVert\paramatom{\vtheta}-\sampling{h}\paramatom{\vtheta}\right\rVert_{L_2(\Bd)}\dd|\mu|(\vtheta)\\
 &\lesssim_{d,s,\kappa}h^{(2s+1)/4}\int_{\Param}\dd{|\mu|}(\vtheta) =h^{(2s+1)/4}\tvnorm{\mu}.
\end{align*}
Taking the infimum over \(\mu\) gives
\begin{align*}
 \lVert u-\sampling{h}u\rVert_{L_2(\Bd)} \lesssim_{d,s,\kappa}h^{(2s+1)/4}\vnorm{u}.
\end{align*}
If \(u|_{X_h}=0\), then Proposition \ref{prop:local-reconstruction} gives \(\sampling{h}u=0\), which yields the second conclusion of the lemma.
This completes the proof.
\end{proof}

\begin{proof}[Proof of Theorem \ref{thm:direct}]
Let \(u=f-\ixhf\in\Vsp\), so that \(u|_{X_h}=0\).
If \(X_h\) satisfies the boundary separation condition with \(\kappa=1/2\), Lemma \ref{lem:operator-error} gives
\begin{align}\label{eq:direct-zero-estimate}
 \lVert u\rVert_{L_2(\Bd)} \lesssim_{d,s}h^{(2s+1)/4}\vnorm{u}.
\end{align}
If \(X_h\) does not satisfy this boundary separation condition, \leavevmode\phantomsection\label{not:contracted-ball}consider the smaller closed ball
\begin{align*}
 K_h:=\{\vx\in\R^d:|\vx|\le1-3h\}.
\end{align*}
For each \(\vz\in K_h\), the definition of the fill distance allows us to choose \(\vy(\vz)\in X_h\) such that \(|\vy(\vz)-\vz|\le h\), \leavevmode\phantomsection\label{not:interior-subset}and we set
\begin{align*}
 Y_h:=\{\vy(\vz):\vz\in K_h\}\subset X_h.
\end{align*}
Since \(X_h\) is finite, so is \(Y_h\), and every \(\vy\in Y_h\) satisfies
\begin{align*}
 |\vy|\le1-2h.
\end{align*}
\leavevmode\phantomsection\label{not:subset-fill-distance}Let the fill distance of \(Y_h\) in \(\Bd\) be denoted by \(\widetilde h\).
For any \(\vx\in\Bd\), if \(|\vx|\le1-3h\), take \(\vz=\vx\); if \(|\vx|>1-3h\), take
\begin{align*}
 \vz=(1-3h)\frac{\vx}{|\vx|}.
\end{align*}
In both cases, \(|\vx-\vz|\le3h\), and consequently
\begin{align*}
 \dist(\vx,Y_h) \le|\vx-\vz|+|\vz-\vy(\vz)| \le4h.
\end{align*}
Thus \(\widetilde h\le4h\), and
\begin{align*}
 |\vy|\le1-2h\le1-\frac12\widetilde h, \qquad \vy\in Y_h.
\end{align*}
This shows that \(Y_h\), relative to its fill distance \(\widetilde h\), satisfies the boundary separation condition with \(\kappa=1/2\).
Since \(Y_h\subset X_h\), we also have \(u|_{Y_h}=0\).
Applying Lemma \ref{lem:operator-error} to the node set \(Y_h\), and using \(\widetilde h\le4h\), again gives \eqref{eq:direct-zero-estimate}.

Finally, \eqref{eq:direct-zero-estimate} and the minimum-norm property of the interpolant imply
\begin{align*}
 \lVert f-\ixhf\rVert_{L_2(\Bd)} &=\lVert u\rVert_{L_2(\Bd)} \lesssim_{d,s}h^{(2s+1)/4}\vnorm{u}\\
 &\lesssim_{d,s}h^{(2s+1)/4} \bigl(\vnorm{f}+\vnorm{\ixhf}\bigr)\\
 &\lesssim_{d,s}h^{(2s+1)/4}\vnorm{f}.
\end{align*}
This completes the proof.
\end{proof}

\begin{proof}[Proof of Theorem \ref{thm:direct-lower}]
Theorem \ref{thm:direct} first yields an upper bound for \(\mathcal R_s(X_h)\).
Indeed, if \(u|_{X_h}=0\), the zero function is a minimum-norm interpolant of the data \(u|_{X_h}\).
Hence, for every such \(u\) with \(\vnorm{u}\le1\), \(\ltnorm{u}\lesssim_{d,s}h^{(2s+1)/4}\), and consequently \(\mathcal R_s(X_h)\lesssim_{d,s}h^{(2s+1)/4}\).
Under the boundary separation condition \eqref{eq:interior-margin}, a matching example can be constructed.
\leavevmode\phantomsection\label{not:sharpness-cap}Fix \(\vxi\in\Sph^{d-1}\), and set \(\delta_h=\kappa h/2\) and \(t_h=1-\delta_h\).
For every \(\vx_i\in X_h\),
\[
 \vxi\cdot\vx_i-t_h \le \abs{\vx_i}-1+\delta_h \le -\frac{\kappa h}{2}<0.
\]
Thus, the normalized ridge atom \(\atom{\vxi}{t_h}\) vanishes on \(X_h\), and \(\vnorm{\atom{\vxi}{t_h}}\le1\).
By Lemma \ref{lem:cap-scale}, \(\norm[\Ltwo]{\atom{\vxi}{t_h}} \asymp_{d,s,\kappa}h^{(2s+1)/4}\), which gives the matching lower bound for \(\mathcal R_s(X_h)\).
\end{proof}

\subsection{Duality and Proofs for the Inverse Theorems}\label{app:inverse-proofs}

\begin{proof}[Proof of Proposition \ref{prop:finite-certificates}]
Write \(S=\dictsynthop{\Lambda}\) and \(G=\wgram{\Lambda}\).
The case of an empty dictionary is immediate.
The given representation has cost \(\vnorm{g}\), and any coefficient vector representing \(g\) has cost at least \(\vnorm{g}\).
Finite-dimensional linear programming duality \citep[Sections~5.2.3--5.2.4]{boydvandenberghe2004convex} therefore gives \eqref{eq:finite-primal-problem}.
The atoms span \(E\), so \(S^*\) is injective on \(E\).
Consequently, the feasible set of the dual problem is closed and bounded in the finite-dimensional space \(E\), hence compact.
The maximum is attained, and \(\mathcal C_\Lambda(g)\) is a nonempty compact convex set.
Strict convexity of the squared Hilbert norm gives a unique minimum-norm certificate.

\leavevmode\phantomsection\label{not:certificate-coordinates}For \(q\in E\), put \(\vect a=S^*q\).
Since \(q=S\vect b\) for some \(\vect b\), we have \(\vect a=G\vect b\in\ran G\).
Conversely, if \(\vect a\in\ran G\), then
\begin{equation*}
 q=SG^\dagger\vect a,\qquad S^*q=GG^\dagger\vect a=\vect a.
\end{equation*}
It is the unique element of \(E\) having the prescribed inner products with the atoms, and
\begin{equation*}
 \ltnorm{q}^2=(G^\dagger\vect a)^{\mathsf T}G(G^\dagger\vect a)=\vect a^{\mathsf T}G^\dagger\vect a.
\end{equation*}
Under \(\norm[\ell^\infty]{\vect a}\le1\), the equality
\begin{equation*}
 \ip{q}{g}_{\Ltwo}=\sum_k c_ka_k=\sum_k|c_k|
\end{equation*}
holds precisely when \(a_k=\operatorname{sgn}(c_k)\) for each nonzero coefficient.
This proves \eqref{eq:certificate-quadratic-program}.

If \(G\) is positive definite, choose \(a_k\) to be the required sign on the support of \(\vc\) and zero elsewhere.
This vector is feasible, and
\begin{equation*}
 \Gamma_\Lambda(g)^2\le\vect a^{\mathsf T}G^{-1}\vect a
 \le\frac{\norm[\ell^2]{\vect a}^2}{\lambda_{\min}(G)}
 \le\frac{\abs{\Lambda}}{\lambda_{\min}(G)}.
\end{equation*}
Taking square roots proves \eqref{eq:certificate-spectral-bound}.
\end{proof}

\begin{proof}[Proof of Theorem \ref{thm:certificate-inverse}]
\leavevmode\phantomsection\label{not:certificate}Take the minimum-norm certificate \(q_j\in\mathcal C_j\), so that \(\ltnorm{q_j}=\Gamma_j\).
Applying \eqref{eq:certificate-increment} with \(g=g_{j+1}\), \(u=g_j\), and \(\Lambda=\jointdict{j}\) and using the triangle inequality gives
\begin{align}\label{eq:minimum-norm-growth}
 m_{j+1}-m_j &\le\ip{q_j}{g_{j+1}-g_j}_{\Ltwo}\le\Gamma_j\norm[\Ltwo]{g_{j+1}-g_j} \le\Gamma_j(e_j+e_{j+1}).
\end{align}
Summing up to \(J-1\) gives
\begin{equation*}
 m_J\le m_1+\sum_{j=1}^{J-1}\Gamma_j(e_j+e_{j+1}).
\end{equation*}
It follows that if \(e_J\to0\) and $\sum_{j=1}^\infty\Gamma_j(e_j+e_{j+1})<\infty$, then \((m_J)\) is bounded, and Proposition~\ref{thm:weak-star-closure} gives \(f\in\Vsp\), with
\begin{equation*}
 \vnorm{f}\le m_1+\sum_{j=1}^\infty\Gamma_j(e_j+e_{j+1}).
\end{equation*}
We now verify these two conditions under the hypotheses of the theorem.
Both error bounds in \eqref{eq:critical-log-rate} imply \(e_j\to0\).

If \(e_j\le Ah_j^\beta\) and \(\beta>\mu\), \eqref{eq:geometric-scales} gives $\Gamma_j(e_j+e_{j+1}) \lesssim_{A,C_\Gamma,\underline\rho,\beta,\mu}h_j^{\beta-\mu}$.
Since \(h_j\le h_1\overline\rho^{\,j-1}\), this is summable.
Under the second error condition, the same argument bounds the summands by \(C(1+j)^{-1-\eps}\), which is also summable.
Finally, \eqref{eq:geometric-scales} gives \(1+j\asymp|\log h_j|\) for sufficiently large \(j\), which yields the stated logarithmic form of the second error condition.
\end{proof}

\begin{proof}[Proof of Remark \ref{cor:approximate-certificates}]
The certificate constraints give
\begin{equation*}
 m_{j+1}-m_j \le\tau_j+\ip{q_j}{g_{j+1}-g_j}_{\Ltwo} \le\tau_j+\ltnorm{q_j}(e_j+e_{j+1}).
\end{equation*}
Summing this inequality and applying Proposition~\ref{thm:weak-star-closure} gives \(f\in\Vsp\) and
\begin{equation*}
 \vnorm{f}\le m_1+\sum_{j=1}^\infty\bigl[\tau_j+\ltnorm{q_j}(e_j+e_{j+1})\bigr],
\end{equation*}
provided that \(e_j\to0\) and
\begin{equation*}
 \sum_{j=1}^\infty\bigl[\tau_j+\ltnorm{q_j}(e_j+e_{j+1})\bigr]<\infty.
\end{equation*}
Under the hypotheses of the remark, \(\sum_j\tau_j<\infty\), and the estimates in the proof of Theorem~\ref{thm:certificate-inverse}, with \(\Gamma_j\) replaced by \(\ltnorm{q_j}\), give \(e_j\to0\) and convergence of the remaining series.
\end{proof}

\begin{proof}[Proof of Corollary \ref{thm:gram-inverse}]
We first explain how a Bernstein estimate for successive differences yields the inverse conclusion.
\leavevmode\phantomsection\label{not:secant-constant}For a given target function \(f\), define the Bernstein constant for differences of interpolants by
\begin{equation*}
 \bernconst{j}{f}:= \sup_{\substack{\minset{j}{f}\ni u\neq v\in\minset{j+1}{f}}} \frac{\vnorm{u-v}}{\norm[\Ltwo]{u-v}}.
\end{equation*}
If no such pair of distinct minimizers exists, we adopt the convention \(\bernconst{j}{f}=0\).
For the chosen sequence \((g_j)\), it suffices to have bounds \(B_j\ge0\) such that
\begin{equation}\label{eq:restricted-bernstein}
 \vnorm{g_{j+1}-g_j}\le B_j\norm[\Ltwo]{g_{j+1}-g_j}.
\end{equation}
One may take \(B_j=\bernconst{j}{f}\) whenever these constants are finite.
By the triangle inequality,
\begin{equation*}
 \vnorm{g_{j+1}-g_j} \le B_j\norm[\Ltwo]{g_{j+1}-g_j} \le B_j(e_j+e_{j+1}).
\end{equation*}
If \(e_j\to0\) and \(\sum_j B_j(e_j+e_{j+1})<\infty\), completeness of \(\Vsp\), established in Appendix~\ref{app:basic-properties}, implies that the series \(\sum_{j=1}^\infty(g_{j+1}-g_j)\) converges in \(\Vsp\), so \((g_j)\) converges in \(\Vsp\) to some \(g\in\Vsp\).
The continuous embedding \(\Vsp\hookrightarrow\Ltwo\) gives \(g_j\to g\) in \(\Ltwo\), whereas \(e_j\to0\) gives \(g_j\to f\) in \(\Ltwo\).
Hence, \(f=g\in\Vsp\).
Letting the partial sums tend to their limit yields
\begin{equation*}
 \vnorm{f} \le\vnorm{g_1}+\sum_{j=1}^\infty B_j(e_j+e_{j+1}).
\end{equation*}
\leavevmode\phantomsection\label{not:bernstein-growth}If \(B_j\le C_Bh_{j+1}^{-\mu}\) for some \(C_B>0\) and \(\mu\ge0\), either error condition in \eqref{eq:critical-log-rate} gives \(e_j\to0\) and convergence of this series, by the estimates in the proof of Theorem~\ref{thm:certificate-inverse}.

For the selected network representations, we obtain suitable bounds \(B_j\) from the weighted Gram matrix.
Consider a finite parameter set \(\Lambda=\{\vtheta_1,\ldots,\vtheta_M\}\).
For any \(v=\dictsynth{\Lambda}{\vc}\), the definition of the weighted atomic cost and the Cauchy--Schwarz inequality give
\begin{equation*}
 \vnorm{v}\le\norm[\ell^1]{\vc} \le\sqrt M\norm[\ell^2]{\vc}.
\end{equation*}
On the other hand, the definition of the weighted Gram matrix and the lower bound for its smallest eigenvalue yield
\begin{equation*}
 \norm[\Ltwo]{v}^2 =\vc^{\mathsf T}\wgram{\Lambda}\vc \ge\mingrameig{\Lambda}\norm[\ell^2]{\vc}^2.
\end{equation*}
Therefore, whenever \(\mingrameig{\Lambda}>0\),
\begin{equation}\label{eq:gram-dictionary-bernstein}
 \vnorm{v} \le \sqrt{\frac{M}{\mingrameig{\Lambda}}}\, \norm[\Ltwo]{v}.
\end{equation}

Combine the representations of \(g_j\) and \(g_{j+1}\) using the atoms indexed by \(\jointdict{j}\), merging the coefficients associated with repeated parameters.
Then \(g_{j+1}-g_j\) belongs to the finite-dimensional space generated by this dictionary.
Applying \eqref{eq:gram-dictionary-bernstein} to \(\jointdict{j}\) and then using the two bounds in \eqref{eq:spectral-condition} gives
\begin{align*}
 \vnorm{g_{j+1}-g_j} &\le \sqrt{\frac{\abs{\jointdict{j}}} {\mingrameig{\jointdict{j}}}} \norm[\Ltwo]{g_{j+1}-g_j}\\
 &\le \sqrt{\frac{C_M}{c_G}}\, h_{j+1}^{-\nu/2-\eta}\norm[\Ltwo]{g_{j+1}-g_j}.
\end{align*}
Thus, \eqref{eq:restricted-bernstein} holds with \(B_j=\sqrt{C_M/c_G}\,h_{j+1}^{-\nu/2-\eta}\).
It remains to verify the summability condition under the assumed error bound.
If \(e_j\le Ah_j^\beta\) and \(\beta>\nu/2+\eta\), then \eqref{eq:geometric-scales} gives
\begin{align*}
 h_{j+1}^{-\nu/2-\eta}(e_j+e_{j+1}) &\le A h_{j+1}^{-\nu/2-\eta}(h_j^\beta+h_{j+1}^\beta)\le A(1+\underline\rho^{-\beta})h_{j+1}^{\beta-\nu/2-\eta}.
\end{align*}
Moreover,
\begin{equation*}
 h_{j+1}^{\beta-\nu/2-\eta} \le h_2^{\beta-\nu/2-\eta} \overline\rho^{\,(j-1)(\beta-\nu/2-\eta)}.
\end{equation*}
Thus, the series converges and \(e_j\to0\), so \(f\in\Vsp\), with
\begin{equation*}
 \vnorm{f}\le\vnorm{g_1}
 +\sqrt{\frac{C_M}{c_G}}\left(1+\underline\rho^{-\beta}\right)
 A\frac{h_2^{\beta-\nu/2-\eta}}{1-\overline\rho^{\,\beta-\nu/2-\eta}}.
\end{equation*}
\end{proof}

We conclude with two one-dimensional examples showing that a Bernstein inequality on \(\netclass{n}\) with a constant depending only on \(n\) cannot hold uniformly over all hidden parameters.
Both examples have fixed width, and their \(L_2\) norms tend to zero while their weighted variation norms tend to infinity.
\leavevmode\phantomsection\label{not:one-dimensional-atom}Write \(\psi_t=\atom{1}{t}\).
By Lemma \ref{lem:weighted-parameter-smoothness}, the parameter map is locally Lipschitz in \(L_2(-1,1)\) when the knot stays away from the boundary.
\leavevmode\phantomsection\label{not:cancellation-example}Therefore, for
\begin{equation*}
 u_\varepsilon :=\varepsilon^{-1/2}(\psi_\varepsilon-\psi_0)
\end{equation*}
we have \(\norm[\Ltwo]{u_\varepsilon}\lesssim_s\varepsilon^{1/2}\to0\).
On the other hand, \(\mathscr D^{s+1}u_\varepsilon\) has Dirac masses of magnitude \(\asymp_s\varepsilon^{-1/2}\) at \(0\) and \(\varepsilon\).
Because the weights at these two knots remain uniformly bounded away from zero, every representation must incur the weighted atomic cost associated with these masses.
Hence, \(\vnorm{u_\varepsilon}\gtrsim_s\varepsilon^{-1/2}\to\infty\).
This example shows how cancellation between two nearby atoms is reflected in the degeneration of the weighted Gram matrix.

\leavevmode\phantomsection\label{not:boundary-example}A second type of degeneration occurs near the boundary.
Specifically, let
\begin{equation*}
 v_\delta:=\delta^{-(2s+1)/8}\psi_{1-\delta}.
\end{equation*}
By Lemma \ref{lem:cap-scale}, \(\norm[\Ltwo]{\psi_{1-\delta}}\asymp_s\delta^{(2s+1)/4}\), and thus \(\norm[\Ltwo]{v_\delta}\asymp_s\delta^{(2s+1)/8}\to0\).
Meanwhile, \(\mathscr D^{s+1}\psi_{1-\delta}\) has mass \(s!\,\delta^{-(2s+1)/4}\) at \(1-\delta\).
Since the weight at this point is exactly \(\delta^{(2s+1)/4}\), it follows that \(\vnorm{v_\delta}\asymp_s\delta^{-(2s+1)/8}\to\infty\).
Thus, a multiple of a single normalized atom approaching the boundary can have a small \(L_2\) norm and a large weighted variation norm, even without cancellation between atoms.
\subsection{A Geometric Sufficient Condition}\label{app:private-patches}
We prove Theorem~\ref{thm:private-inverse} using the local hyperplane separation condition \eqref{eq:private-patch}.
For a finite collection of distinct hyperplanes intersecting \(\Bd\), the condition holds for some positive radius: on each hyperplane, choose an interior point outside the other hyperplanes, and then take sufficiently small disjoint balls.
The theorem requires a quantitative lower bound on these radii as the fill distance tends to zero.
The condition is sufficient but not necessary for a positive lower spectral bound.
For example, the one-dimensional atoms \((x)_+^s\) and \((-x)_+^s\) have a positive definite Gram matrix, but their associated hyperplanes coincide.
Moreover, pairwise disjointness of the balls implies \(\abs{\jointdict{j}}\privrad{j}^{\,d}\lesssim_d1\).
If \(\abs{\jointdict{j}}\asymp h_{j+1}^{-\nu}\), then necessarily \(\kappa_0\ge\nu/d\).

\begin{proof}[Proof of Theorem \ref{thm:private-inverse}]
We first establish a Bernstein estimate for the dictionary indexed by a finite parameter set \(\Lambda=\{\vtheta_1,\ldots,\vtheta_M\}\) satisfying the local hyperplane separation condition with radius \(\rho\).
Consider any
\begin{equation*}
 v=\sum_{k=1}^M a_k(\vxi_k\cdot\vx-t_k)_+^s, \qquad c_k=\wstar{t_k}a_k,
\end{equation*}
\leavevmode\phantomsection\label{not:private-ball}and let \(B_k=B(\vy_k,\rho)\).
Fix \(k\).
By the local hyperplane separation condition, for each \(\ell\ne k\), the continuous function \(\vxi_\ell\cdot\vx-t_\ell\) does not vanish on the connected set \(B_k\), so its sign is constant on \(B_k\).
Therefore, the atom \((\vxi_\ell\cdot\vx-t_\ell)_+^s\) is either identically zero on \(B_k\) or is a polynomial of degree at most \(s\) there.
Since \(\vy_k\in H_k\), on \(B_k\) we may write
\begin{equation}\label{eq:private-local-decomposition}
 v(\vx)=P_k(\vx)+a_k\bigl(\vxi_k\cdot(\vx-\vy_k)\bigr)_+^s, \qquad P_k\in\Polys,
\end{equation}
with \(\Polys\) as in Appendix~\ref{app:direct-proof}.

We next estimate the coefficient \(a_k\) in \eqref{eq:private-local-decomposition}.
\leavevmode\phantomsection\label{not:polynomial-distance}Define, on the unit ball \(\Bd\),
\begin{equation*}
 \delta_{d,s}:=\inf_{P\in\Polys} \norm[\Ltwo]{(x_1)_+^s-P}.
\end{equation*}
Since \(\Polys\) is a finite-dimensional closed subspace of \(L_2(\Bd)\) and \((x_1)_+^s\notin\Polys\), we have \(\delta_{d,s}>0\).
Indeed, if some \(P\in\Polys\) were equal to \((x_1)_+^s\) almost everywhere, then \(P\) would vanish almost everywhere on the open set \(\{\vx\in\Bd:x_1<0\}\), and hence would be the zero polynomial.
This contradicts its equality to \(x_1^s\) almost everywhere on \(\{\vx\in\Bd:x_1>0\}\).
\leavevmode\phantomsection\label{not:private-rotation}Choose an orthogonal matrix \(\matr{Q}_k\) such that \(\matr{Q}_k\vect{e}_1=\vxi_k\), and make the change of variables \(\vx=\vy_k+\rho\matr{Q}_k\vz\).
Restricting \eqref{eq:private-local-decomposition} to \(B_k\), scaling to the unit ball, and then taking the distance to \(\Polys\), we obtain
\begin{align*}
 \rho^{-d/2}\norm[L_2(B_k)]{v} &=\norm[L_2(\Bd)]{v(\vy_k+\rho\matr{Q}_k\,\cdot)} \ge \abs{a_k}\rho^s\delta_{d,s}.
\end{align*}
Consequently,
\begin{equation}\label{eq:private-coefficient-bound}
 \abs{a_k} \lesssim_{d,s}\rho^{-s-d/2}\norm[L_2(B_k)]{v}.
\end{equation}
Because the balls are pairwise disjoint and contained in \(\Bd\), summing the squares of \eqref{eq:private-coefficient-bound} gives
\begin{equation*}
 \norm[\ell^2]{\vect a} \lesssim_{d,s}\rho^{-s-d/2}\norm[\Ltwo]{v}.
\end{equation*}
For \(-1\le t_k\le1\), we have \(0\le\wstar{t_k}\le2^{(d+2s)/4}\).
Therefore, \(c_k=\wstar{t_k}a_k\) implies
\begin{equation*}
 \norm[\ell^2]{\vc} \lesssim_{d,s}\rho^{-s-d/2}\norm[\Ltwo]{v}.
\end{equation*}
On the other hand, \(v=\sum_{k=1}^M c_k\paramatom{\vtheta_k}\), and hence
\begin{equation}\label{eq:private-dictionary-bernstein}
 \vnorm{v}\le\norm[\ell^1]{\vc} \le\sqrt M\norm[\ell^2]{\vc} \lesssim_{d,s}\sqrt M\,\rho^{-s-d/2}\norm[\Ltwo]{v}.
\end{equation}
Notice that \(\vy_k\in H_k\cap\Bd\) implies \(t_k<1\), so \(\wstar{t_k}>0\).
Finally, for any \(\vc\in\R^M\), let \(a_k=c_k/\wstar{t_k}\).
Applying \eqref{eq:private-coefficient-bound} to \(\dictsynth{\Lambda}{\vc}\) and using the pairwise disjointness of the balls gives
\begin{equation*}
 \vc^{\mathsf T}\wgram{\Lambda}\vc =\norm[\Ltwo]{\dictsynth{\Lambda}{\vc}}^2 \gtrsim_{d,s}\rho^{2s+d}\norm[\ell^2]{\vc}^2.
\end{equation*}
The Rayleigh quotient characterization therefore yields \(\mingrameig{\Lambda}\gtrsim_{d,s}\rho^{2s+d}\).

Now combine the finite network representations of \(g_j\) and \(g_{j+1}\) on \(\jointdict{j}\).
Applying \eqref{eq:private-dictionary-bernstein} to \(\jointdict{j}\) and using the bounds on the dictionary cardinality and the radius of the balls gives
\begin{align*}
 \vnorm{g_{j+1}-g_j} &\lesssim_{d,s} \abs{\jointdict{j}}^{1/2}\privrad{j}^{-s-d/2} \norm[\Ltwo]{g_{j+1}-g_j}\\
 &\lesssim_{d,s,C_M,c_\rho} h_{j+1}^{-\nu/2-\kappa_0(s+d/2)} \norm[\Ltwo]{g_{j+1}-g_j}.
\end{align*}
By the summability argument in the proof of Corollary~\ref{thm:gram-inverse}, it follows that \(f\in\Vsp\) whenever \(e_j\to0\) and
\begin{equation*}
 \sum_{j=1}^\infty h_{j+1}^{-\nu/2-\kappa_0(s+d/2)}(e_j+e_{j+1})<\infty.
\end{equation*}
If \(e_j\le Ah_j^\beta\) and \(\beta>\nu/2+\kappa_0(s+d/2)\), then \eqref{eq:geometric-scales} gives
\begin{align*}
 &h_{j+1}^{-\nu/2-\kappa_0(s+d/2)}(e_j+e_{j+1}) \qquad\lesssim_{A,\underline\rho,\beta} h_{j+1}^{\beta-\nu/2-\kappa_0(s+d/2)}.
\end{align*}
The terms on the right are bounded by a geometric sequence with ratio \(\overline\rho^{\,\beta-\nu/2-\kappa_0(s+d/2)}<1\).
The assumed error bound also gives \(e_j\to0\), so the preceding summability argument completes the proof.
\end{proof}
\clearpage
\subsection{Table of Notation}\label{app:notation}

The following table lists the principal notation.
Page numbers refer to the corresponding definitions or introductions in the text.

\begingroup \footnotesize \setlength{\tabcolsep}{4pt} \setlength{\LTleft}{0pt plus 1fill} \setlength{\LTright}{0pt plus 1fill} \renewcommand{\arraystretch}{1.10} \renewcommand{\thetable}{\thesubsection.\arabic{table}}
\begin{longtable}{@{}>{\raggedright\arraybackslash}p{.21\textwidth} >{\raggedright\arraybackslash}p{.20\textwidth} >{\raggedright\arraybackslash}p{\dimexpr.59\textwidth-6\tabcolsep-0.95cm\relax} >{\centering\arraybackslash}p{0.95cm}@{}}
\caption{Principal notation.}\label{tab:notation}\\
\toprule
\textbf{Notation} & \textbf{Name} & \textbf{Meaning} & \textbf{Page}\\
\midrule
\endfirsthead
\toprule
\textbf{Notation} & \textbf{Name} & \textbf{Meaning} & \textbf{Page}\\
\midrule
\endhead
\bottomrule
\endlastfoot
\(a=O(b)\) & Asymptotic bound & \(\abs{a}\le C\abs{b}\) eventually in the stated limit, with \(C>0\) independent of the limiting parameter. & \pageref{tab:notation}\\[2pt]
\(a\lesssim_{d,s}b\), \(a\gtrsim_{d,s}b\) & Comparison bounds & For \(a,b\ge0\), \(a\le C_{d,s}b\) or \(b\le C_{d,s}a\), respectively, where \(C_{d,s}>0\) depends only on \(d,s\). & \pageref{tab:notation}\\[2pt]
\(a\asymp_{d,s}b\) & Comparability & \(a\lesssim_{d,s}b\) and \(b\lesssim_{d,s}a\), possibly with different constants. & \pageref{tab:notation}\\[2pt]
\(\#A\) & Cardinality & Number of elements in the set \(A\). & \pageref{tab:notation}\\[2pt]
\(\Bd\) & Unit ball & \(\{\vx\in\R^d:\abs{\vx}<1\}\). & \pageref{not:domain}\\[2pt]
\(L_p(\Omega)\), \(\norm[L_p(\Omega)]{\cdot}\) & Lebesgue space and norm & \(\norm[L_p(\Omega)]{f}=(\int_\Omega\abs{f}^p\dd{\vx})^{1/p}\) for \(1\le p<\infty\); the essential supremum for \(p=\infty\). & \pageref{not:lp-spaces}\\[2pt]
\(H^\beta(\Omega)\) & Sobolev space & The \(L_2\) Sobolev space of order \(\beta\). & \pageref{not:sobolev-spaces}\\[2pt]
\(\Sph^{d-1}\), \(\Param\) & Unit sphere and parameter domain & \(\Sph^{d-1}=\{\vxi\in\R^d:\abs{\vxi}=1\}\), \(\Param=\Sph^{d-1}\times[-1,1]\). & \pageref{not:parameter-domain}\\[2pt]
\(\Msp\), \(\tvnorm{\cdot}\) & Radon measures and total variation & Finite signed Radon measures on \(\Param\), with \(\tvnorm{\mu}=\abs{\mu}(\Param)\). & \pageref{not:radon-measures}\\[2pt]
\(\norm[\ell^1]{\vc}\) & Coefficient norm & \(\sum_k\abs{c_k}\). & \pageref{eq:finite-primal-problem}\\[2pt]
\(\wstar{t}\) & Boundary weight & \((1-t)^{(d+2s)/4}\). & \pageref{not:weight}\\[2pt]
\(\atom{\vxi}{t}\), \(\paramatom{\vtheta}\) & Normalized atom & \((\vxi\cdot\vx-t)_+^s/\wstar{t}\) for \(t<1\), and zero at \(t=1\). & \pageref{eq:normalized-cap-atom}\\[2pt]
\(\synthop\) & Synthesis operator & \(\synth{\mu}=\int_{\Param}\paramatom{\vtheta}\dd{\mu}(\vtheta)\). & \pageref{not:synthesis}\\[2pt]
\(\Vsp\), \(\vnorm{\cdot}\) & Weighted variation space and norm & \(\Vsp=\ran\synthop\), with \(\vnorm{f}=\inf_{\synth{\mu}=f}\tvnorm{\mu}\). & \pageref{not:weighted-space}\\[2pt]
\(\mathfrak M_X(f)\) & Minimizer set & Minimum-norm solutions subject to \(g|_X=f|_X\); abbreviated as \(\minset{j}{f}\) for \(X=X_j\). & \pageref{not:minimizer-set}\\[2pt]
\(\interp{X}{f}\) & Minimum-norm interpolant & An arbitrary element of \(\mathfrak M_X(f)\). & \pageref{not:interpolant}\\[2pt]
\(h,\ h_j,\ q_X\) & Fill and separation distances & \(h\) is the supremum of the distance to the nearest node; \(h_j=h_{X_j,\Bd}\), and \(q_X\) is half the minimum distance between distinct nodes. & \pageref{eq:fill-distance}\\[2pt]
\(\mathcal R_s(X_h)\) & Radius for zero data & The supremum of \(\ltnorm{u}\) over \(\vnorm{u}\le1\) and \(u|_{X_h}=0\). & \pageref{not:sampling-radius}\\[2pt]
\(g_j,e_j,m_j\) & Interpolants, errors, and norms & \(g_j\in\minset{j}{f}\), \(e_j=\norm[\Ltwo]{f-g_j}\), and \(m_j=\vnorm{g_j}\). & \pageref{not:inverse-errors}\\[2pt]
\(\dict{j},\jointdict{j}\) & Parameter sets & \(\dict{j}\) indexes a representation of \(g_j\) attaining the minimum weighted atomic cost; \(\jointdict{j}=\dict{j}\cup\dict{j+1}\). & \pageref{eq:minimum-cost-representations}\\[2pt]
\(E_j\) & Span of the combined dictionary & \(\spann\{\paramatom{\vtheta}:\vtheta\in\jointdict{j}\}\). & \pageref{not:joint-dictionary}\\[2pt]
\(\wgram{\Lambda}\) & Weighted Gram matrix & The matrix with entries \(\ip{\paramatom{\vtheta_k}}{\paramatom{\vtheta_\ell}}_{\Ltwo}\). & \pageref{not:dictionary-matrices}\\[2pt]
\(\mathcal C_j,\Gamma_j\) & Dual certificates and minimum norm & \(\mathcal C_j\) is the certificate set in \eqref{eq:finite-certificate-set}; \(\Gamma_j=\min_{q\in\mathcal C_j}\ltnorm{q}\). & \pageref{eq:finite-certificate-set}\\[2pt]
\(\netclass{n}\) & Network class & Shallow \(\relu^s\) networks of width at most \(n\). & \pageref{not:network-class}\\[2pt]
\(\bestnerr{f}\) & Best approximation error & \(\inf_{g\in\netclass{n}}\ltnorm{f-g}\). & \pageref{not:network-error}\\
\end{longtable}
\endgroup
\clearpage

\section*{Data availability}
No data was used for the research described in the article.

\FloatBarrier
\bibliographystyle{plain}
\bibliography{references}
\end{document}